\documentclass[12pt]{article}

\usepackage{amsmath}
\usepackage{amssymb}
\usepackage[dvips]{graphicx}

\begin{document}
\title{{\bf 
{\large Finite slope cyclic surgeries along toroidal Brunnian links
and generalized Properties P and R}}
\footnotetext[0]{%
2010 {\it Mathematics Subject Classification}: 
11R27, 57M25, 57M27, 57Q10, 57R65.\\
{\it Keywords}:
Milnor link; Dehn surgery; Property P; Property R; 
Reidemeister torsion; Alexander polynomial.
}}
\author{{\footnotesize 
Teruhisa KADOKAMI}}
\date{{\footnotesize April 16, 2015}}
\maketitle

\newcommand{\circlenum}[1]{{\ooalign{%
\hfill$\scriptstyle#1$\hfill\crcr$\bigcirc$}}}

\newcommand{\svline}[1]{\multicolumn{1}{|c}{#1}}
\newfont{\bg}{cmr10 scaled\magstep4}
\newcommand{\bigzerol}{\smash{\hbox{\bg 0}}}
\newcommand{\bigzerou}{\smash{\lower1.7ex\hbox{\bg 0}}}

\newcommand{\bsquare}{\hbox{\rule{6pt}{6pt}}}
\newcommand{\qed}{\hbox{\rule[-2pt]{3pt}{6pt}}}
\newcommand{\Int}{\mathrm{Int}\ \! }
\newcommand{\Ker}{\mathrm{Ker}\ \! }
\newcommand{\Ig}{\mathrm{Im}\ \! }
\newcommand{\aug}{\mathrm{aug}\ \!}
\newcommand{\pj}{\mathrm{pr}\ \!}
\newcommand{\Tor}{\mathrm{Tor}\ \!}
\newcommand{\Spin}{\mathrm{Spin}\ \!}
\newcommand{\Eul}{\mathrm{Eul}\ \!}
\newcommand{\Vect}{\mathrm{Vect}\ \!}
\newcommand{\HULL}{\mathrm{HULL}\ \!}
\newcommand{\real}{\mathrm{real}\ \!}
\newcommand{\rank}{\mathrm{rank}\ \!}
\newcommand{\ord}{\mathrm{ord}\ \!}
\newcommand{\Sign}{\mathrm{Sign}\ \!}
\newcommand{\Hom}{\mathrm{Hom}\ \!}
\newcommand{\ad}{\mathrm{ad}\ \!}
\newcommand{\Det}{\mathrm{Det}\ \!}
\newcommand{\lk}{\mathrm{lk}\ \!}
\newcommand{\pt}{\mathrm{pt}}
\newcommand{\al}{$\alpha$}
\newcommand{\dis}{\displaystyle}

\newtheorem{df}{Definition}[section]
\newtheorem{lm}[df]{Lemma}
\newtheorem{theo}[df]{Theorem}
\newtheorem{re}[df]{Remark}
\newtheorem{pr}[df]{Proposition}
\newtheorem{ex}[df]{Example}
\newtheorem{co}[df]{Corollary}
\newtheorem{cl}[df]{Claim}
\newtheorem{qu}[df]{Question}
\newtheorem{pb}[df]{Problem}
\newtheorem{cj}[df]{Conjecture}

\makeatletter
\renewcommand{\theequation}{%
\thesection.\arabic{equation}}
\@addtoreset{equation}{section}
\makeatother

\begin{abstract}
{\footnotesize 
\setlength{\baselineskip}{10pt}
\setlength{\oddsidemargin}{0.25in}
\setlength{\evensidemargin}{0.25in}
Let $M_{\lambda}$ be the $\lambda$-component Milnor link.
For $\lambda \ge 3$,
we determine completely when a finite slope surgery
along $M_{\lambda}$ yields 
a lens space including $S^3$ and $S^1\times S^2$,
where {\it finite slope surgery} implies that
a surgery coefficient of every component is not $\infty$.
For $\lambda =3$ (i.e.\ the Borromean rings), 
there are three infinite sequences of finite slope surgeries
yielding lens spaces.
For $\lambda \ge 4$,
any finite slope surgery does not yield a lens space.
As a corollary, 
$M_{\lambda}$ for $\lambda \ge 3$
does not yield both $S^3$ and $S^1\times S^2$
by any finite slope surgery.
We generalize the results for the cases of 
{\it Brunnian type links} and toroidal Brunnian type links
(i.e.\ Brunnian type links including essential tori in the link complement).
Our main tools are Alexander polynomials and Reidemeister torsions.
Moreover we characterized 
toroidal Brunnian links and toroidal Brunnian type links in some senses.
}
\end{abstract}

\tableofcontents

\section{Introduction}~\label{sec:intro}
A $\lambda$-component link $L$ with $\lambda \ge 1$ is a {\it Brunnian link}
if every proper sublink of $L$ is a trivial link.
For $\lambda \ge 3$, $L$ is algebraically split
(i.e.\ the linking number of every $2$-component sublink of $L$ is $0$),
and for $\lambda \ge 1$, every component of $L$ is an unknot.
We assume that for $\lambda=1$, $L$ is the unknot, and
for $\lambda=2$, $L$ is algebraically split.

\medskip

Let
$M_{\lambda}
=K_1\cup \cdots \cup K_{\lambda}$ be the $\lambda$-component Milnor link
\cite{Mi1, Mi2} (Figure 1).
In particular, $M_1$ is the unknot, $M_2$ is the Hopf link, 
and $M_3$ is the Borromean rings.
Since Dehn surgeries along the unknot and the Hopf link are well-known
(for example, see \cite{KS}), we assume $\lambda \ge 3$ for $M_{\lambda}$.
As properties of $M_{\lambda}$,
it is a Brunnian link, 
it can be constructed from the Hopf link by a sequence of Bing doubles
(cf.\ Subsection \ref{ssec:geogen}),
the components $K_1$ and $K_2$, and $K_{\lambda-1}$ and $K_{\lambda}$
are interchangeable, and the order of the components
$(K_1, \ldots, K_{\lambda})$ can be reversed into
$(K_{\lambda}, \ldots, K_1)$.
The Alexander polynomial of $M_{\lambda}$ is as follows
(see also Lemma \ref{lm:BrunnianAlex} (3) and Remark \ref{re:MAlex}) :
\begin{equation}\label{eq:Alex}
\left\{
\begin{array}{cccl}
{\mit \Delta}_{M_3}(t_1, t_2, t_3) & \doteq & (t_1-1)(t_2-1)(t_3-1) & 
\quad (\lambda=3),
\medskip\\
{\mit \Delta}_{M_{\lambda}}(t_1, \ldots, t_{\lambda}) & = & 0 &
\quad (\lambda \ge 4).
\end{array}
\right.
\end{equation}

\begin{figure}[htbp]
\begin{center}
\includegraphics[scale=0.7]{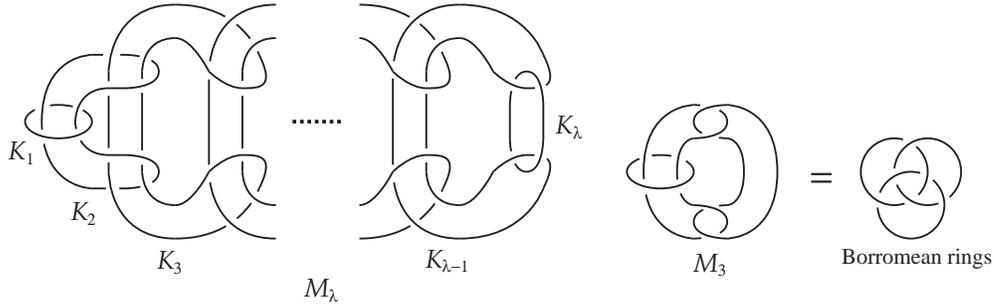}
\label{M}
\caption{$\lambda$-component Milnor link $M_{\lambda}$}
\end{center}
\end{figure}

Let $L=K_1\cup \cdots \cup K_{\lambda}$ be a $\lambda$-component link.
We denote the result of
$(r_1, \ldots, r_{\lambda})$-surgery along $L$ by
$(L; r_1, \ldots, r_{\lambda})$ where
$r_i\in \mathbb{Q}\cup \{\infty, \emptyset\}\ (i=1, \ldots, \lambda)$
is a surgery coefficient of $K_i$.
A Dehn surgery producing a $3$-manifold with cyclic fundamental group 
is called a {\it cyclic surgery}.
In particular, a Dehn surgery producing a lens space is called a {\it lens surgery}.
We define that an $(r_1, \ldots, r_{\lambda})$-surgery along $L$
is a {\it finite slope surgery} if $r_i\ne \infty, \emptyset$ 
for every $i$\ $(i=1, \ldots, \lambda)$.
When we set $r_i=p_i/q_i\in \mathbb{Q}\cup \{\infty\}$,
we assume $p_i, q_i\in \mathbb{Z}$, $\gcd(p_i, q_i)=1$ and
$\pm 1/0=\infty$.
The lens space of type $(p, q)$
is the result of $p/q$-surgery along the unknot,
which is denoted by $L(p, q)$.
In particular, $L(1, q)\cong S^3$ and $L(0, \pm 1)\cong S^1\times S^2$.
We obtain the following theorems:

\begin{theo}\label{th:MT1}
Let $Y=(M_3; p_1/q_1, p_2/q_2, p_3/q_3)$ be
the result of a finite slope surgery along $M_3$.
Then $Y$ is a lens space if and only if
one of the following holds:
\begin{enumerate}
\item[(1)]
$(p_1/q_1, p_2/q_2)
=(\varepsilon, \varepsilon)$ and $|\varepsilon p_3-6q_3|=1,$ or

\item[(2)]
$(p_1/q_1, p_2/q_2)
=(\varepsilon, 2\varepsilon)$ and $|\varepsilon p_3-4q_3|=1,$ or

\item[(3)]
$(p_1/q_1, p_2/q_2)
=(\varepsilon, 3\varepsilon)$ and $|\varepsilon p_3-3q_3|=1,$

\end{enumerate}
and the cases that indices of $(p_i, q_i)\ (i=1, 2, 3)$ are permuted arbitrarily,
where $\varepsilon=1$ or $-1$.
Moreover, if (1), then $Y=L(p_3, 4\varepsilon q_3)$, 
if (2), then $Y=L(2p_3, \varepsilon (8q_3-p_3))$, and 
if (3), then $Y=L(3p_3, \varepsilon (3q_3-2p_3))$.
\end{theo}

\begin{theo}\label{th:MT2}
For $\lambda \ge 4$, 
any finite slope surgery along $M_{\lambda}$ does not yield a lens space.
\end{theo}

For a non-trivial knot $K$ in $S^3$,
$K$ has {\it Property P} if the fundamental group of 
the result of any finite slope surgery 
($1/q$-surgery with $q\ne 0$) along $K$ is a non-trivial group, 
and $K$ has {\it Property R} if the result of $0$-surgery 
along $K$ is not homeomorphic to $S^1\times S^2$.
Note that by L.~Moser \cite{Mos} 
which determines the results of Dehn surgery along torus knots,
and by the cyclic surgery theorem due to
M.~Culler, M.~Gordon, J.~Luecke and P.~Shalen \cite{CGLS},
$q$ can be restricted to $q=\pm 1$
in the statement of Property P.
We define that a non-trivial knot $K$ in $S^3$ has 
{\it Property P$'$} if the result of any finite slope surgery along $K$ 
is not $S^3$, which is a stronger property than Property P.
Since the Poincar\'e conjecture has been
settled affirmatively by G.~Perelman \cite{Pe1, Pe2, Pe3},
Property P and Property P$'$ are equivalent.
However, as Property P has been introduced to prove
the Poincar\'e conjecture, we cannot use the statement of Property P$'$
when we try to give an alternative proof of the Poincar\'e conjecture.

\medskip

The {\it Property P conjecture} 
({\it Property R conjecture}, respectively)
for a knot is that a non-trivial knot in $S^3$ has Property P (Property R, respectively).
As for the Property P conjecture for a knot, C.~McA.~Gordon and J.~Luecke \cite{GL} 
gave a great partial affirmative answer, and
P.~Kronheimer and T.~Mrowka \cite{KM} proved affirmatively for every non-trivial knot.
As for the Property R conjecture for a knot,
D.~Gabai \cite{Gb2} proved affirmatively.

\medskip

We define their link versions.
For a $\lambda$-component non-trivial link $L$ in $S^3$,
$L$ has {\it Property P} ({\it Property FP}, respectively)
if the fundamental group of the result of any Dehn surgery 
(finite slope surgery, respectively) along $L$ 
whose surgery coefficients are not all $\infty$
(i.e.\ $(r_1, \ldots, r_{\lambda})\ne (\infty, \ldots, \infty)$)
is non-trivial, 
$L$ has {\it Property P$'$} ({\it Property FP$'$}, respectively)
if the result of any Dehn surgery 
(finite slope surgery, respectively) along $L$ 
whose surgery coefficients are not all $\infty$ is not $S^3$, 
and $L$ has {\it Property R} ({\it Property FR}, respectively)
if the result of any Dehn surgery (finite slope surgery, respectively)
along $L$ is not homeomorphic to $S^1\times S^2$.
It is easy to see that not all non-trivial link
has Property P, Property FP, Property P$'$, Property FP$'$, 
Property R and Property FR
even if we assume that the link is non-split.
For example, the Hopf link is a such link.
By the same reason as above,
Property P and Property P$'$, and
Property FP and Property FP$'$ are equivalent, but we should 
distinguish them for an alternative proof of the Poincar\'e conjecture.

\medskip

A link $L$ in $S^3$ has a {\it unique exterior}
if there exists a link $L'$ and a homeomorphism $h$ between 
the complements of $L$ and $L'$, then
$h$ preserves the peripheral structures of $L$ and $L'$
(i.e.\ $L$ and $L'$ are the same links).
By C.~McA.~Gordon and J.~Luecke \cite{GL}, 
every knot in $S^3$ has this property.
It is easy to see that a link with a non-split unknotted component
does not have a unique exterior.
In particular, every Milnor link does not have a unique exterior.
Recently A.~Kawauchi \cite{Kw3} proved 
by using the imitation theory \cite{Kw1, Kw2}
that the statement of the Poincar\'e conjecture
is equivalent to the statement that
a link $L$ in $S^3$ having a unique exterior has Property FP.

\medskip

From Theorem \ref{th:MT1} and Theorem \ref{th:MT2}, 
we have a corollary
which implies that the Milnor links with $\lambda \ge 3$
has Property FP$'$ and Property FR.

\begin{co}\label{co:propPR}
For $\lambda \ge 3$, 
any finite slope surgery along $M_{\lambda}$ does not yield 
both $S^3$ and $S^1\times S^2$.
\end{co}

Since the proof of Corollary \ref{co:propPR} comes down to
the proof of non-triviality of a certain knot in $S^3$
by the result of P.~Kronheimer and T.~Mrowka \cite{KM},
$M_{\lambda}$ with $\lambda \ge 3$ has Property FP.
B.~Mangum and T.~Stanford \cite{MS} proved that
every non-trivial Brunnian link is determined by its complement.
We notice that it does not imply that every non-trivial Brunnian link has a unique exterior,
but implies that if a link $L$ in $S^3$ has a homeomorphic complement to a Brunnian link $L'$,
then $L$ is equivalent to $L'$ or $L$ is not a Brunnian link.
Therefore Corollary \ref{co:propPR} does not contribute to 
an alternative proof of the Poincar\'e conjecture.
However it provides the cyclic surgery problem of links
with sensitive examples.

\medskip

A $\lambda$-component link in $S^3$ is a {\it Brunnian type link} 
if the link is an algebraically split link such that every component is an unknot
and the Alexander polynomial of every proper sublink with at least two components is zero
(Subsection \ref{ssec:alggen}).
In Section \ref{sec:gen}, we generalize our results for the cases of 
Brunnian type links 
(Theorem \ref{th:MT3}, Theorem \ref{th:MT4} and Corollary \ref{co:const}), 
and toroidal Brunnian type links
(Corollary \ref{co:cond} (2))
which are Brunnian type links including essential tori in the link complements
(Subsection \ref{ssec:geogen}).
Moreover we characterized 
toroidal Brunnian links and toroidal Brunnian type links in some senses.

\medskip

In Section \ref{sec:R-tor}, we define the Reidemeister torsions, and
state their surgery formulae and properties of Alexander polynomials
which are basic tools in the present paper.
In Section \ref{sec:key}, we show some key lemmas 
for the proofs of Theorem \ref{th:MT1}, Theorem \ref{th:MT2}, and
their generalizations in Section \ref{sec:gen} 
by the techniques in \cite{Kd1, Kd2, Kd3, KMS, Tr1, Tr2} on the Reidemeister torsions.
In Section \ref{sec:3-comp}, we prove Theorem \ref{th:MT1}.
In Section \ref{sec:many-comp}, we prove Theorem \ref{th:MT2} and Corollary \ref{co:propPR}.
A necessary condition for existence of a lens surgery is obtained from 
a lemma in Section \ref{sec:key}, 
and non-existence of a lens surgery is shown by
results due to Y.~Ni \cite{Ni} (Theorem \ref{th:Ni}) 
on fiberedness of a knot yielding a lens space, and
D.~Gabai \cite{Gb1} (Theorem \ref{th:Gabai})
on minimal genus Seifert surfaces for a plumbed link, and
existence of an essential torus in the surgered manifold.
Corollary \ref{co:propPR} is shown by
the proofs of Theorem \ref{th:MT1} and Theorem \ref{th:MT2},
and results due to 
P.~Kronheimer and T.~Mrowka \cite{KM} on Property P for a knot, 
D.~Gabai \cite{Gb1} on Property R for a knot.
We remark that we cannot prove Theorem \ref{th:MT2} 
by only using the Reidemeister torsions and
the Casson-Walker invariants.
This implies that an invariant deduced from the knot Floer homology
due to L.~Nicolaescu \cite{Nc} cannot prove it.
However the knot Floer homology itself can prove it
by Ni's theorem above.

\section{Reidemeister torsion}~\label{sec:R-tor}
We define the Reidemeister torsions, and
state their surgery formulae and properties of Alexander polynomials
which are basic tools in the present paper.

\subsection{Surgery formulae}~\label{ssec:surgery}
Let $X$ be a finite CW complex with $H=H_1(X)$,
$R$ an integral domain, and
$\psi : \mathbb{Z}[H]\to Q(R)$ a ring homomorphism
where $Q(R)$ is the quotient field of $R$.
Then $\tau^{\psi}(X)\in Q(R)$ is 
the {\it Reidemeister torsion} of $X$ related with $\psi$
which is determined up to multiplication of $\pm \psi(h)\ (h\in H)$
(cf. \cite{Tr1, Tr2}).
We do not give a precise definition here.
For $A$ and $B\in Q(R)$,
we denote $A\doteq B$ if $A=\pm \psi(h) B$ for some $h\in H$.
If $\psi$ is the identity map, then we denote $\tau(X):=\tau^{\psi}(X)$ for simplicity.

\medskip

The Reidemeister torsion of $X$ is defined from 
a cell chain complex $\mathbf{C}_{\ast}$ of the maximal abelian covering 
of $X$ over a ring $\mathbb{Z}[H]$.
We can take a basis $\mathbf{e}$ of $\mathbf{C}_{\ast}$
as an ordered oriented lifts of cells of $X$.
Then the value of $\tau^{\psi}(X)$ is determined
as a unique element in $Q(R)$.
A base change with the determinant $1$ preserves
the value of $\tau^{\psi}(X)$, and
a difference of two bases is described by an element
of $\pm H_1(X)$.
We define the basis $\mathbf{e}$ of $\mathbf{C}_{\ast}$
up to base changes with the determinant $1$
a {\it combinatorial Euler structure} of $X$.
The concept is equivalent to
a {\it spin$^c$-structure} if $X$ is a compact oriented $3$-manifold.
A great advantage of the concept is to unify separated values
$\{ \tau^{\psi}(X)\}$ to an element of $Q(\mathbb{Z}[H])$
via the Chinese Remainder Theorem (Suntzu Theorem), and
a troublesome ambiguity vanishes.
The unified value is called the {\it Reidemeister-Turaev torsion}
or the {\it Reidemeister torsion with a combinatorial Euler structure}
or the {\it maximal abelian torsion} (cf. \cite{Tr2}).

\begin{pr}\label{pr:torEx}
\begin{enumerate}
\item[(1)]
We set $H_1(S^1)=\langle t\rangle \cong \mathbb{Z}$.
Then we have $\tau(S^1)\doteq (t-1)^{-1}$.

\item[(2)]
We have $\tau(S^1\times S^1)\doteq 1$.

\item[(3)]
Let $L=K_1\cup \cdots \cup K_{\lambda}$ be a $\lambda$-component link
in an integral homology $3$-sphere, and 
$E_L$ the complement of $L$.
We set $H_1(E_L)=\langle t_1, \ldots, t_{\lambda}\rangle 
\cong \mathbb{Z}^{\lambda}$.
Then we have
$$
\tau(E_L)\doteq 
\left\{
\begin{array}{cl}
{\displaystyle
\frac{{\mit \Delta}_L(t_1)}{t_1-1}} & (\lambda=1),\medskip\\
{\mit \Delta}_L(t_1, \ldots, t_{\lambda}) & (\lambda \ge 2).
\end{array}
\right.$$
\end{enumerate}
\end{pr}

The following is coming from the excision property of the homology groups.

\begin{lm}\label{lm:surgery1}
{\rm (excision)}\ 
Let $M_1$ and $M_2$ be connected compact $3$-manifolds
with non-empty tori boundaries, and
$M=M_1\cup M_2$ a union of $M_1$ and $M_2$
such that $T=M_1\cap M_2=\partial M_1\cap \partial M_2$ is a torus.
Let $i_1 : M_1\to M$, $i_2 : M_2\to M$, 
$j_1 : T\to M_1$ and $j_2 : T\to M_2$ be natural inclusions.
Let $F$ be a field,
$\psi : \mathbb{Z}[H_1(M)]\to F$ a ring homomorphism,
and $\psi_1=\psi \circ (i_1)_{\ast}$, $\psi_2=\psi \circ (i_2)_{\ast}$
and $\rho=\psi \circ (i_1)_{\ast}\circ (j_1)_{\ast}=\psi \circ (i_2)_{\ast}\circ (j_2)_{\ast}$
the induced homomorphisms.
Suppose that $\rho$ is not a trivial map.
Then we have
$$\tau^{\psi}(M)\doteq
\tau^{\psi_1}(M_1)\cdot \tau^{\psi_2}(M_2).$$
\end{lm}

Let $H$ be a finite cyclic group with order $p\ge 2$,
and $t$ a generator.
For a divisor $d\ge 2$ of $p$, 
let $\zeta_d$ be a primitive $d$-th root of unity,
$\psi_d : \mathbb{Z}[H]\to \mathbb{Q}(\zeta_d)$ a ring homomorphism
such that $\psi_d(t)=\zeta_d$.
From Lemma \ref{lm:surgery1}, we have 
the following surgery formulae
which are suitable for studying the results of Dehn surgeries
along links in an integral homology $3$-sphere.

\begin{lm}\label{lm:surgery2}
{\rm (surgery formula)}\ 
\begin{enumerate}
\item[(1)]
Let
$L=K_1\cup \cdots \cup K_{\lambda}$ be a $\lambda$-component link
with $\lambda \ge 2$
in an integral homology $3$-sphere,
$\ell_i'\ (i=1, \ldots, \lambda)$ a core of attaching solid torus, and
$Y$ the result of Dehn surgery along $L$ with $\sharp (H_1(Y))\ge 2$.
Let $F$ be a field, and
$\psi : \mathbb{Z}[H_1(Y)]\to F$ a ring homomorphism.
If $\psi \left( [\ell_i']\right) \ne 1$ for every $i=1, \ldots, \lambda$, 
then we have
$$\tau^{\psi}(Y) \doteq 
{\mit \Delta}_L(\psi([m_1]), \ldots, \psi([m_{\lambda}]))
\cdot \prod_{i=1}^{\lambda}\left( \psi([\ell_i'])-1\right)^{-1}
\in F.$$

\item[(2)]
Let $K$ be a knot in an integral homology $3$-sphere,
$Y=(K; p/q)\ (p\ge 2)$ the result of $p/q$-surgery along $K$,
and $t$ a generator of
$H_1(Y)\cong \mathbb{Z}/p\mathbb{Z}$
corresponding to a meridian of $K$.
Then we have
$$\tau^{\psi_d}(Y)\doteq
{\mit \Delta}_K(\zeta_d)
(\zeta_d-1)^{-1}(\zeta_d^{{\bar q}}-1)^{-1}$$
where $q{\bar q}\equiv 1\ (\mathrm{mod}\ \! p).$
\end{enumerate}
\end{lm}

For example, we have
\begin{equation}\label{eq:lens}
\tau^{\psi_d}(L(p, q))\doteq
(\zeta_d-1)^{-1}(\zeta_d^{{\bar q}}-1)^{-1}.
\end{equation}

The Torres formula for the Alexander polynomial
is a kind of surgery formulae.
Let $L=K_1\cup \cdots \cup K_{\lambda} \cup K_{\lambda+1}\ 
(\lambda \ge 1)$ be an oriented $(\lambda+1)$-component link
in a homology $3$-sphere,
$L'=K_1\cup \cdots \cup K_{\lambda}$ a $\lambda$-component sublink,
and $k_i=\mathrm{lk}\ \! (K_i, K_{\lambda+1})\ (i=1, \ldots, \lambda)$.

\begin{lm}\label{lm:Torres}{\rm (Torres formula \cite{To})}
$${\it \Delta}_L(t_1, \ldots, t_{\lambda}, 1)\doteq
\left\{
\begin{array}{cl}
{\displaystyle \frac{t_1^{k_1}-1}{t_1-1}{\mit \Delta}_{L'}(t_1)} & (\lambda=1),
\medskip\\
(t_1^{k_1}\cdots t_{\lambda}^{k_{\lambda}}-1)
{\mit \Delta}_{L'}(t_1, \ldots, t_{\lambda}) & (\lambda \ge 2).
\end{array}
\right.$$
\end{lm}

\begin{lm}\label{lm:dual}{\rm (duality \cite{Mi3, Tr1})}
Let $L=K_1\cup \cdots \cup K_{\lambda}$ be a $\lambda$-compnent link
in a homology $3$-sphere.
Then we have
$${\mit \Delta}_L(t_1)=t_1^e{\mit \Delta}_L(t_1^{-1})\quad
(\lambda=1)$$
for some even integer $e$, and
$${\mit \Delta}_L(t_1, \ldots, t_{\lambda})=
(-1)^{\lambda}t_1^{e_1}\cdots t_{\lambda}^{e_{\lambda}}
{\mit \Delta}_L(t_1^{-1}, \ldots, t_{\lambda}^{-1})
\quad (\lambda \ge 2)$$
where $e_i\equiv 1+\sum_{j\ne i}\mathrm{lk}\ \! (K_i, K_j)\ 
(\mathrm{mod}\ \! 2)$\ $(i=1, \ldots, \lambda)$.
\end{lm}

\subsection{$d$-norm}~\label{ssec:norm}
About algebraic fields,
we refer the reader to L.~C.~Washington \cite{Was} for example.
For an element $x$ in the $d$-th cyclotomic field $\mathbb{Q}(\zeta_d)$,
the {\it $d$-norm} of $x$ is defined as 
$$N_d(x)=\prod_{\sigma \in \mathrm{Gal}\ \! 
(\mathbb{Q}(\zeta_d)/\mathbb{Q})}
\sigma(x),$$
where $\mathrm{Gal}\ \! (\mathbb{Q}(\zeta_d)/\mathbb{Q})$
is the Galois group related with a Galois extension
$\mathbb{Q}(\zeta_d)$ over $\mathbb{Q}$.
The order of $\mathrm{Gal}\ \! (\mathbb{Q}(\zeta_d)/\mathbb{Q})$
is the Euler function $\varphi : \mathbb{N}\to \mathbb{N}$ at $d$.
In particular, $\varphi(1)=\varphi(2)=1$, and $\varphi(d)$ is even for $d\ge 3$.
The following is well-known.

\begin{pr}~\label{pr:norm}
\begin{enumerate}
\item[(1)]
If $x\in \mathbb{Q}(\zeta_d)$, then $N_d(x)\in \mathbb{Q}$.
The map $N_d : \mathbb{Q}(\zeta_d)\setminus \{0\}\to 
\mathbb{Q}\setminus \{0\}$ is a group homomorphism.

\item[(2)]
If $x\in \mathbb{Z}[\zeta_d]$, then $N_d(x)\in \mathbb{Z}$.
\end{enumerate}
\end{pr}

By easy calculations, we have the following.

\begin{lm}~\label{lm:cyclotomic}
\begin{enumerate}
\item[(1)]
${\displaystyle
N_d(\pm \zeta_d)=
\left\{
\begin{array}{cl}
\pm 1 & (d=2),\\
1 & (d\ge 3).
\end{array}
\right.
}$

\item[(2)]
${\displaystyle
N_d(1-\zeta_d)=
\left\{
\begin{array}{cl}
\ell & (\mbox{$d$ is a power of a prime $\ell \ge 2$}),\\
1 & (\mbox{otherwise}).
\end{array}
\right.
}$
\end{enumerate}
\end{lm}

%



\section{Key Lemmas}~\label{sec:key}
We show some key lemmas 
for the proofs of Theorem \ref{th:MT1}, Theorem \ref{th:MT2}, and
their generalizations in Section \ref{sec:gen} 
by the techniques in \cite{Kd1, Kd2, Kd3, KMS, Tr1, Tr2}
on the Reidemeister torsions.

\subsection{Alexander polynomial of algebraically split links}~\label{ssec:asplit}
Let $L=K_1\cup \cdots \cup K_{\lambda}$ be 
an oriented $\lambda$-component algebraically split link with $\lambda \ge 2$.
By Lemma \ref{lm:Torres}, we have
\begin{equation}\label{eq:Alex0}
{\mit \Delta}_L(t_1, \ldots, t_{\lambda})
=(t_1-1)\cdots (t_{\lambda}-1)f(t_1, \ldots, t_{\lambda})
\end{equation}
where $f(t_1, \ldots, t_{\lambda})\in 
\mathbb{Z}[t_1^{\pm 1}, \ldots, t_{\lambda}^{\pm 1}]$.
By Lemma \ref{lm:dual}, we may assume that
\begin{equation}\label{eq:dual1}
f(t_1, \ldots, t_{\lambda})=f(t_1^{-1}, \ldots, t_{\lambda}^{-1}).
\end{equation}
We add one component $K$ to $L$
such that $L_i=K_i\cup K$\ $(i=1, \ldots, \lambda)$
is the connected sum of $K_i$ and the Hopf link,
where the linking number of $K_i$ and $K$ is $1$.
Then we have
\begin{equation}\label{eq:Alex1}
{\mit \Delta}_{L_i}(t_i, t)
\doteq {\mit \Delta}_{K_i}(t_i).
\end{equation}
We set $\overline{L}=L\cup K$.
By Lemma \ref{lm:Torres} and (\ref{eq:Alex0}), we have
\begin{equation}\label{eq:Torres1}
\begin{matrix}
{\mit \Delta}_{\overline{L}}(t_1, \ldots, t_{\lambda}, t)
& = & (t_1\cdots t_{\lambda}-1)(t_1-1)\cdots (t_{\lambda}-1)
f(t_1, \ldots, t_{\lambda})\medskip\\
& & +(t-1)g(t_1, \ldots, t_{\lambda}, t)\hfill
\end{matrix}
\end{equation}
where $g(t_1, \ldots, t_{\lambda}, t)\in 
\mathbb{Z}[t_1^{\pm 1}, \ldots, t_{\lambda}^{\pm 1}, t^{\pm 1}]$, 
and we may replace $f(t_1, \ldots, t_{\lambda})$ and $g(t_1, \ldots, t_{\lambda}, t)$
with $-f(t_1, \ldots, t_{\lambda})$ and $-g(t_1, \ldots, t_{\lambda}, t)$.
By Lemma \ref{lm:dual}, we may assume that
\begin{equation}\label{eq:dual2}
{\mit \Delta}_{\overline{L}}(t_1, \ldots, t_{\lambda}, t)
=(-1)^{\lambda+1}t_1^2\cdots t_{\lambda}^2t^a
{\mit \Delta}_{\overline{L}}(t_1^{-1}, \ldots, t_{\lambda}^{-1}, t^{-1})
\end{equation}
where $a\equiv \lambda+1\ (\mathrm{mod}\ \! 2)$.
We set $I_{\lambda}=\{1, 2, \ldots, \lambda\}$.
For $I=\{i_1, i_2, \ldots, i_{\mu}\}\subset I_{\lambda}$,
we set $L_I=K_{i_1}\cup \cdots \cup K_{i_s}$,
$\overline{L}_I=L_I\cup K$, and
$g_I\in \mathbb{Z}[t_{i_1}^{\pm 1}, \ldots, t_{i_{\mu}}^{\pm 1}, t^{\pm 1}]$ 
is obtained by substituting 
$t_j=1$ for all $j\in I_{\lambda}\setminus I$
to $g(t_1, \ldots, t_{\lambda}, t)$.

\medskip

By (\ref{eq:Alex1}) and (\ref{eq:dual2}),
if $I=\{i\}$\ $(s=1)$, then we may take
\begin{equation}\label{eq:Torres2}
{\mit \Delta}_{\overline{L}_I}(t_i, t)
={\mit \Delta}_{K_i}(t_i)
\end{equation}
where ${\mit \Delta}_{K_i}(t_i)=t_i^2{\mit \Delta}_{K_i}(t_i^{-1})$
and ${\mit \Delta}_{K_i}(1)=1$.
If $2\le s \le \lambda$\ $(\lambda \ge 2)$, then we may take
\begin{equation}\label{eq:Torres3}
\begin{matrix}
{\mit \Delta}_{\overline{L}_I}(t_{i_1}, \ldots, t_{i_s}, t)
& = & {\displaystyle \left( \prod_{i\in I}t_i-1\right)\prod_{i\in I}(t_i-1)
f_I(t_{i_1}, \ldots, t_{i_s})}\hfill \medskip\\
&  & {\displaystyle +(t-1)
g_I'(t_{i_1}, \ldots, t_{i_s}, t)}\hfill
\end{matrix}
\end{equation}
where
\begin{equation*}
{\mit \Delta}_{L_I}(t_{i_1}, \ldots, t_{i_s})
=\prod_{i\in I}(t_i-1)f_I(t_{i_1}, \ldots, t_{i_s}),
\end{equation*}
\begin{equation*}
f_I(t_{i_1}, \ldots, t_{i_s})=f_I(t_{i_1}^{-1}, \ldots, t_{i_s}^{-1})
\in \mathbb{Z}[t_{i_1}^{\pm 1}, \ldots, t_{i_{\mu}}^{\pm 1}]
\end{equation*}
and
\begin{equation*}
g_I'(t_{i_1}, \ldots, t_{i_s}, t)
\in \mathbb{Z}[t_{i_1}^{\pm 1}, \ldots, t_{i_{\mu}}^{\pm 1}, t^{\pm 1}].
\end{equation*}
We set $f_I=f_I(t_{i_1}, \ldots, t_{i_s})$ and
$g_I'=g_I'(t_{i_1}, \ldots, t_{i_s}, t)$.
We remark that 
$f_{I_{\lambda}}=f=f(t_1, \ldots, t_{\lambda})$ and
$g_{I_{\lambda}}'=g=g(t_1, \ldots, t_{\lambda}, t)$.

\begin{lm}\label{lm:Torres2}
Under the situation above, for $1\le s \le \lambda-1$, we have
\begin{equation*}
g_I=(t-1)^{\lambda-s-1}
{\mit \Delta}_{\overline{L}_I}(t_{i_1}, \ldots, t_{i_s}, t).
\end{equation*}
\end{lm}

\noindent
{\it Proof.}\ 
By applying Lemma \ref{lm:Torres} repeatedly, we have the result.
\qed

\bigskip

We exapand the $g(t_1, \ldots, t_{\lambda}, t)$-part
in (\ref{eq:Torres1}) with normalization as follows:

\begin{lm}\label{lm:form}
For $\lambda \ge 2$, we have
\begin{equation*}
\begin{matrix}
g(t_1, \ldots, t_{\lambda}, t)
& = & {\displaystyle (t-1)^{\lambda-2}
\prod_{i=1}^{\lambda}{\mit \Delta}_{K_i}(t_i)}\hfill \medskip\\
& & {\displaystyle +\sum_{
{\scriptstyle I\subset I_{\lambda}}
\atop
{\scriptstyle 2\le s=|I| \le \lambda-1}}
\prod_{i\in I}(t_i-1)
(t-1)^{\lambda-s-1}}\hfill\medskip\\
& & {\displaystyle 
\cdot \left\{\left( \prod_{i\in I}t_i-1\right)f_I
+(t-1)h_I\right\}}\hfill\medskip\\
& & {\displaystyle +\prod_{i=1}^{\lambda}(t_i-1)h}\hfill
\end{matrix}
\end{equation*}
where $|I|$ is the number of elements of $I$, and
$f_I\in \mathbb{Z}[t_{i_1}^{\pm 1}, \ldots, t_{i_{\mu}}^{\pm 1}]$,
$h_I\in \mathbb{Z}[t_{i_1}^{\pm 1}, \ldots, t_{i_{\mu}}^{\pm 1}, t^{\pm 1}]$
and $h\in \mathbb{Z}[t_1^{\pm 1}, \ldots, t_{\lambda}^{\pm 1}, t^{\pm 1}]$.

\end{lm}

For the proof of Lemma \ref{lm:form},
we refer the reader to \cite[Lemma 4.2]{Kd3}.

\begin{co}\label{co:g}
$g(1, \ldots, 1, t)=(t-1)^{\lambda-2}$.
\end{co}

\begin{re}\label{re:g}
{\rm
By Corollary \ref{co:g}, the normalization of $g(t_1, \ldots, t_{\lambda}, t)$
is uniquely determined, and 
the normalization of $f(t_1, \ldots, t_{\lambda})$
is also uniquely determined.
}
\end{re}

\subsection{Key Lemmas}~\label{ssec:key}
Let $L$ and $\overline{L}$ be the same links as in Subsection \ref{ssec:asplit}.

\medskip

We set $Y=(L; p_1/q_1, \ldots, p_{\lambda}/q_{\lambda})$
and $\overline{Y}=(\overline{L}; p_1/q_1, \ldots, p_{\lambda}/q_{\lambda}, \emptyset)$
with $q_i\ne 0\ (i=1, \cdots, \lambda)$.
We calculate a presentation matrix of 
the first homologies $H_1(Y)$ and $H_1(\overline{Y})$.
Let $m_i$ and $l_i$ be a meridian and a longitude of $K_i$ ($i=1, \ldots, \lambda$),
and $m$ and $l$ a meridian and a longitude of $K$.
Let $T_i$ be an attaching solid torus of $K_i$, 
$T$ an attaching solid torus of $K$, 
$m_i'$ and $l_i'$ a meridian and a core of $T_i$,
and $m'$ and $l'$ a meridian and a core of $T$.
We denote the representing element of a loop $\gamma$ in the first homology by $[\gamma]$.
Then we have
\begin{eqnarray*}
H_1(E_L) & = & \langle [m_1], \ldots, [m_{\lambda}]\rangle 
\cong \mathbb{Z}^{\lambda}\medskip\\
H_1(E_{\overline{L}}) & = & \langle [m_1], \ldots, [m_{\lambda}], [m]\rangle 
\cong \mathbb{Z}^{\lambda+1}
\end{eqnarray*}
and
\begin{eqnarray*}
H_1(Y) & = & 
\langle [m_1], [l_1]\ldots, [m_{\lambda}], [l_{\lambda}]\ |\ 
[m_i]^{p_i}[l_i]^{q_i}=1, [l_i]=1\ (i=1, \cdots, \lambda)\rangle \medskip\\
 & \cong & \langle [m_1], \ldots, [m_{\lambda}]\ |\ [m_i]^{p_i}=1\ (i=1, \cdots, \lambda)\rangle
\medskip\\
 & \cong & {\displaystyle \bigoplus_{i=1}^{\lambda}\mathbb{Z}/p_i\mathbb{Z}}.
\end{eqnarray*}
Hence $H_1(Y)\cong \mathbb{Z}/p\mathbb{Z}$
if and only if $p=|p_1\cdots p_{\lambda}|$ and 
$\gcd(p_i, p_j)=1\ (1\le i\ne j\le \lambda)$.
From now on, we assume $p\ge 2$.

\medskip

\begin{equation}\label{eq:Ybar}
\begin{array}{ccl}
H_1(\overline{Y}) & = & 
\left\langle 
\begin{matrix}
[m_1], [l_1], \ldots, [m_{\lambda}], [l_{\lambda}], \hfill\\
[m_1'], [l_1']\ldots, [m_{\lambda}'], [l_{\lambda}'],\hfill\\
[m], [l] \hfill
\end{matrix}
\ \left|\
\begin{matrix}
[m_i']=[m_i]^{p_i}[l_i]^{q_i}=1, \hfill\\
[l_i']=[m_i]^{r_i}[l_i]^{s_i}\ (p_is_i-q_ir_i=-1), \hfill\\
[l_i]=[m]\ (i=1, \cdots, \lambda),\hfill\\
[l]=[m_1]\cdots [m_{\lambda}] \hfill
\end{matrix}
\right.
\right\rangle
\medskip\\
& \cong & 
\langle 
[m_1], \ldots, [m_{\lambda}], [m]\ |\ 
[m_i]^{p_i}[m]^{q_i}=1\ (i=1, \cdots, \lambda)
\rangle
\medskip\\
& \cong & \mathbb{Z}.
\end{array}
\end{equation}

We set $t_i=[m_i]$, $t=[m]$ in $H_1(Y)$ and $H_1(\overline{Y})$, and
a generator of $H_1(\overline{Y})$ as $T$.
Then a group homomorphism
$\rho :  H_1(\overline{Y})\to \langle T\rangle \cong \mathbb{Z}$
defined by
\begin{equation}\label{eq:rho1}
\rho(t_i)=T^{\frac{q_ip}{p_i}}, \quad \rho(t)=T^{-p}
\end{equation}
is a well-defined isomorphism.
Then we have
\begin{equation}\label{eq:rho2}
\rho([l_i'])=T^{\frac{p}{p_i}}, \quad 
\rho([l])=\rho(t_1\cdots t_{\lambda})
=T^{\frac{q_1p}{p_1}+\cdots +\frac{q_{\lambda}p}{p_{\lambda}}}.
\end{equation}

\medskip

We note that the result of $\infty$-surgery along $K$ from $\overline{Y}$ is $Y$, and
since it corresponds to $T^{-p}=1$ in (\ref{eq:rho1}), 
$$H_1(Y)\cong \langle T\ |\ T^p=1\rangle \cong \mathbb{Z}/p\mathbb{Z}$$
is recovered.

\medskip

For fixed $k$, we set $f_k=f_k(t_k)=f(1, \ldots, 1, t_k, 1, \ldots, 1)$
which is obtained by substituting $t_j=1\ (j\ne k)$ to $f(t_1, \ldots, t_{\lambda})$.
Then we have the following:

\begin{lm}\label{lm:key1}
Suppose that $L$ is a $\lambda$-component Brunnian link with $\lambda \ge 3$.
For a divisor $d\ge 2$ of $p_k$, 
we have
$$\tau^{\psi_d}(Y)\doteq
\left\{ f_k(\zeta_d)
\left(
\prod_{\scriptsize{j\ne k}
\atop
\scriptsize{1\le j\le \lambda}}
q_j
\right)
(\zeta_d-1)^2
+(-1)^{\lambda}
\left(
\prod_{\scriptsize{j\ne k}
\atop
\scriptsize{1\le j\le \lambda}}
p_j
\right)
\zeta_d\right\}
\cdot (\zeta_d-1)^{-1}(\zeta_d^{{\bar q}_k}-1)^{-1}$$
where $q_k{\bar q}_k\equiv 1\ (\mathrm{mod}\ \! p_k)$.
\end{lm}

\noindent
{\it Proof.}\ 
Since $L$ is a Brunnian link, 
we have ${\mit \Delta}_{K_i}(t_i)=t_i$ and $f_I=0$ in Lemma \ref{lm:form}.
By Lemma \ref{lm:surgery2} (1), Lemma \ref{lm:form}, (\ref{eq:rho1}) and (\ref{eq:rho2}), 
we have
\begin{equation}\label{eq:tauYbar1}
\begin{matrix}
\tau(\overline{Y})
& \doteq & {\displaystyle
{\mit \Delta}_{\overline{L}}
(\rho(t_1), \ldots, \rho(t_{\lambda}), \rho(t))
\prod_{i=1}^{\lambda}(\rho([l_i'])-1)^{-1}}\hfill
\medskip\\
& \doteq &{\displaystyle
\left(T^{\frac{q_1p}{p_1}+\cdots +\frac{q_{\lambda}p}{p_{\lambda}}}-1\right)
\prod_{i=1}^{\lambda}
\frac{T^{\frac{q_ip}{p_i}}-1}{T^{\frac{p}{p_i}}-1}\cdot
f\left(
T^{\frac{q_1p}{p_1}}, \ldots, T^{\frac{q_{\lambda}p}{p_{\lambda}}}
\right)
}\hfill
\medskip\\
&  & {\displaystyle
+(-1)^{\lambda}T^{p\lambda+\frac{q_1p}{p_1}+\cdots +\frac{q_{\lambda}p}{p_{\lambda}}}
\prod_{\scriptsize{j\ne k}
\atop
\scriptsize{1\le j\le \lambda}}
\frac{T^p-1}{T^{\frac{p}{p_j}}-1}\cdot
\left(T^{\frac{p}{p_k}}-1\right)^{-1}
}\hfill
\medskip\\
&  & {\displaystyle
+(T^{-p}-1)\sum_{
{\scriptstyle I\subset I_{\lambda}}
\atop
{\scriptstyle 2\le s=|I| \le \lambda-1}}
\prod_{i\in I}\frac{T^{\frac{q_ip}{p_i}}-1}{T^{\frac{p}{p_i}}-1}
\prod_{j\in I_{\lambda}\setminus I}
\frac{T^{-p}-1}{T^{\frac{p}{p_j}}-1}
}\hfill
\medskip\\
&  & {\displaystyle
\cdot h_I\left(
T^{\frac{q_{i_1}p}{p_{i_1}}}, \ldots, T^{\frac{q_{i_{\mu}}p}{p_{i_{\mu}}}}, T^{-p}
\right)
}\hfill
\medskip\\
&  & {\displaystyle
+(T^{-p}-1)\prod_{i=1}^{\lambda}
\frac{T^{\frac{q_ip}{p_i}}-1}{T^{\frac{p}{p_i}}-1}\cdot
h\left(
T^{\frac{q_1p}{p_1}}, \ldots, T^{\frac{q_{\lambda}p}{p_{\lambda}}}, T^{-p}
\right).
}\hfill
\medskip\\
\end{matrix}
\end{equation}

\medskip

Since $p/p_k$ and $q_k$ are relatively prime to $p_k$,
$\frac{q_kp}{p_k}$ is also relatively prime to $p_k$.
We define a ring homomorphism
$\overline{\psi}_d : \mathbb{Z}[T, T^{-1}]\to \mathbb{Q}(\zeta_d)$ by
$$\overline{\psi}_d(T)=\zeta_d^{\overline{\frac{q_kp}{p_k}}}$$
where 
$\frac{q_kp}{p_k}\cdot \overline{\frac{q_kp}{p_k}}\equiv 1\ (\mathrm{mod}\ \! p_k)$.
Then we have
$$\overline{\psi}_d\left( T^{\frac{q_ip}{p_i}}\right)
=\overline{\psi}_d\left( T^{\frac{p}{p_i}}\right)=1\ (i\ne k),\ 
\overline{\psi}_d\left( T^{\frac{q_kp}{p_k}}\right)=\zeta_d,\ 
\overline{\psi}_d\left( T^{\frac{p}{p_k}}\right)=\zeta_d^{{\bar q}_k},\ 
\overline{\psi}_d\left( T^p\right)=1,$$
$$
\begin{array}{ccccll}
{\displaystyle
\overline{\psi}_d\left( \frac{T^{\frac{q_ip}{p_i}}-1}{T^{\frac{p}{p_i}}-1}\right)}
& = & {\displaystyle
\mathrm{sign}(q_i)\cdot
\overline{\psi}_d\left( 
\sum_{j=1}^{|q_i|-1}T^{\frac{jp}{p_i}}\right)} & = & q_i & (i\ne k),
\medskip\\
{\displaystyle
\overline{\psi}_d\left( \frac{T^{p}-1}{T^{\frac{p}{p_i}}-1}\right)}
& = & {\displaystyle
\mathrm{sign}(p_i)\cdot
\overline{\psi}_d\left( 
\sum_{j=1}^{|p_i|-1}T^{\frac{jp}{p_i}}\right)} & = & p_i & (i\ne k),
\end{array}$$
and by (\ref{eq:tauYbar1}), we have
\begin{equation}\label{eq:tauYbar2}
\tau^{\overline{\psi}_d}(\overline{Y})\doteq
\left\{ f_k(\zeta_d)
\left(
\prod_{\scriptsize{j\ne k}
\atop
\scriptsize{1\le j\le \lambda}}
q_j
\right)
(\zeta_d-1)^2
+(-1)^{\lambda}
\left(
\prod_{\scriptsize{j\ne k}
\atop
\scriptsize{1\le j\le \lambda}}
p_j
\right)
\zeta_d\right\}
\cdot (\zeta_d^{{\bar q}_k}-1)^{-1}.
\end{equation}
By Lemma \ref{lm:surgery1}, 
$\psi_d([l])=\zeta_d$, and (\ref{eq:tauYbar2}),
we have
$$\begin{array}{ccl}
\tau^{\psi_d}(Y) & \doteq &
\tau^{\overline{\psi}_d}(\overline{Y})\cdot (\zeta_d-1)^{-1}
\medskip\\
& = &
{\displaystyle
\left\{ f_k(\zeta_d)
\left(
\prod_{\scriptsize{j\ne k}
\atop
\scriptsize{1\le j\le \lambda}}
q_j
\right)
(\zeta_d-1)^2
+(-1)^{\lambda}
\left(
\prod_{\scriptsize{j\ne k}
\atop
\scriptsize{1\le j\le \lambda}}
p_j
\right)
\zeta_d\right\}}
\bigskip\\
& &
\cdot (\zeta_d-1)^{-1}(\zeta_d^{{\bar q}_k}-1)^{-1}.
\end{array}$$
\qed

\begin{lm}\label{lm:key2}
In the same situation as in Lemma \ref{lm:key1}, we suppose $Y$ is a lens space.
Then we have the following:

\begin{enumerate}
\item[(1)]
If $f_k$ is constant with $f_k\ge 1$ and $|p_k|\ge 5$, then we have
$f_k=1$, $|q_j|=1$\ $(j\ne k)$,  
and
${\displaystyle
\prod_{\scriptsize{j\ne k}
\atop
\scriptsize{1\le j\le \lambda}}
p_j=\eta}$, $2\eta$ or $3\eta$
where ${\displaystyle
\eta=\prod_{\scriptsize{j\ne k}
\atop
\scriptsize{1\le j\le \lambda}}
q_j=1}$ or $-1$.

\item[(2)]
If $f_k=0$, then we have $|p_j|=1$\ $(j\ne k)$.

\end{enumerate}
\end{lm}

\noindent
{\it Proof.}\ 
Since $Y$ is a lens space and (\ref{eq:lens}), we have
$$\tau^{\psi_d}(Y)\doteq (\zeta_d^a-1)^{-1}(\zeta_d^b-1)^{-1}$$
for some $a$ and $b$,
where $\gcd(a, p_k)=\gcd(b, p_k)=1$.
Then by Lemma \ref{lm:key1}, we have an equation
\begin{equation}\label{eq:lenseq}
\left\{ f_k(\zeta_d)
\left(
\prod_{\scriptsize{j\ne k}
\atop
\scriptsize{1\le j\le \lambda}}
q_j
\right)
(\zeta_d-1)^2
+(-1)^{\lambda}
\left(
\prod_{\scriptsize{j\ne k}
\atop
\scriptsize{1\le j\le \lambda}}
p_j
\right)
\zeta_d\right\}\cdot
(\zeta_d^a-1)(\zeta_d^b-1)\doteq (\zeta_d-1)(\zeta_d^c-1)
\end{equation}
where $c\equiv \pm {\bar q}_k\ (\mathrm{mod}\ \! p_k)$.

\medskip

By fixing a combinatorial Euler structure,
we can lift the equation above
to the equation in $\mathbb{Z}[u, u^{-1}]/(u^{|p_k|-1}+\cdots +u+1)$
as follows : 
\begin{equation}\label{eq:lift1}
\left\{ f_k(u)
\left(
\prod_{\scriptsize{j\ne k}
\atop
\scriptsize{1\le j\le \lambda}}
q_j
\right)
(u-1)^2
+(-1)^{\lambda}
\left(
\prod_{\scriptsize{j\ne k}
\atop
\scriptsize{1\le j\le \lambda}}
p_j
\right)
u\right\}\cdot
(u^a-1)(u^b-1)=\pm u^l (u-1)(u^c-1)
\end{equation}
and it deduces the following : 
\begin{equation}\label{eq:lift2}
\left\{ f_k(u)
\left(
\prod_{\scriptsize{j\ne k}
\atop
\scriptsize{1\le j\le \lambda}}
q_j
\right)
(u-1)^2
+(-1)^{\lambda}
\left(
\prod_{\scriptsize{j\ne k}
\atop
\scriptsize{1\le j\le \lambda}}
p_j
\right)
u\right\}\cdot
\frac{u^a-1}{u-1}\cdot
\frac{u^b-1}{u-1}
=\pm u^l \cdot \frac{u^c-1}{u-1}.
\end{equation}

\medskip

We set the part $\{\cdots\}$ in (\ref{eq:lift1}) and (\ref{eq:lift2}) as $F(u)$.
If $f_k$ is constant, then $F(u)$ is of degree $2$.
We set $F(u)=mu^2+nu+m$ ($m, n\in \mathbb{Z}$).
By taking the $d$-norm (cf.\ Subsection \ref{ssec:norm}) of both sides of (\ref{eq:lenseq}),
we have $(m, n)\ne (0, 0)$ and $\gcd(m, n)=1$
(The $d$-norm of the part $\{\cdots\}$ in (\ref{eq:lenseq}) is $1$).

\medskip

\noindent
(1) Suppose that $f_k$ is constant with $f_k\ge 1$ and $|p_k|\ge 5$.

\medskip

\noindent
(i) The case that $p_k$ is even.

\medskip

We can take $1\le a, b, c\le |p_k|/2-1$ are odd, 
the degree of the lefthand side of (\ref{eq:lift2}) is $\le |p_k|-2$, and
the degree of the righthand side of (\ref{eq:lift2}) is $\le |p_k|/2-2$.
Since both lefthand side and righthand side of (\ref{eq:lift2})
are symmetric polynomials, they are identical as elements in $\mathbb{Z}[u, u^{-1}]$
up to multiplications of $\pm u^l$.
Hence $F(u)$ is a divisor of $(u^c-1)/(u-1)$, and
$F(u)$ is a cyclotomic polynomial of degree $2$
$$\pm (u^2+u+1),\quad
\pm (u^2+1)\ \mbox{or}\ 
\pm (u^2-u+1).$$
Therefore we have the result.

\medskip

\noindent
(ii) The case that $p_k$ is odd.

\medskip

We can take $1\le a, b, c\le |p_k|-1$ satisfying
$2\le a+b\le |p_k|-1$ and $a+b\equiv c+1\ (\mathrm{mod}\ \! 2)$.
If $a+b\le |p_k|-2$, then by the same argument as (i), we have the result.
Suppose that $a+b=|p_k|-1$ 
(i.e.\ $a=b=(|p_k|-1)/2$ and $c$ is odd with $1\le c\le |p_k|-2$).
We calculate the lefthand side of (\ref{eq:lift1}), and deform it modulo
$(u^{|p_k|-1}+\cdots +u+1)$.
\begin{equation}\label{eq:deform}
\begin{array}{l}
(mu^2+nu+m)\left( u^{\frac{|p_k|-1}{2}}-1\right)^2= 
(mu^2+nu+m)\left(u^{|p_k|-1}-2u^{\frac{|p_k|-1}{2}}+1
\right)
\medskip\\
=mu^{|p_k|+1}+nu^{|p_k|}+mu^{|p_k|-1}
-2mu^{\frac{|p_k|+3}{2}}-2nu^{\frac{|p_k|+1}{2}}-2mu^{\frac{|p_k|-1}{2}}
\medskip\\
+mu^2+nu+m
\medskip\\
=u^{\frac{|p_k|+1}{2}}\left\{
-2n-2m(u+u^{-1})+m\left( u^{\frac{|p_k|-3}{2}}+u^{-\frac{|p_k|-3}{2}}\right)
+(m+n)\left( u^{\frac{|p_k|-1}{2}}+u^{-\frac{|p_k|-1}{2}}\right)
\right\}
\medskip\\
=u^{\frac{|p_k|+1}{2}}\left\{
(-m-3n)+(-3m-n)(u+u^{-1})
-(m+n)(u^2+u^{-2})-
\cdots
\right.
\medskip\\
\left.
-(m+n)\left( u^{\frac{|p_k|-5}{2}}+u^{-\frac{|p_k|-5}{2}}\right)
-n\left( u^{\frac{|p_k|-3}{2}}+u^{-\frac{|p_k|-3}{2}}\right)
\right\}.
\end{array}
\end{equation}

\medskip

\noindent
(a) The case $3\le c\le |p_k|-4$ ($|p_k|\ge 7$).

\medskip

The lefthand side of (\ref{eq:lift1}) is $\pm u^l(u^{c+1}-u^c-u+1)$.
By (\ref{eq:deform}), we have $-m-3n=0$, and
$(m, n)=(3, -1)$ or $(-3, 1)$.
However $-3m-n=-8$ or $8$.
Hence we have no root for $(m, n)$.

\medskip

\noindent
(b) The case $c=1$.

\medskip

The lefthand side of (\ref{eq:lift1}) is $\pm u^l(u^2-2u+1)$.
By (\ref{eq:deform}), we have $-m-3n=\mp 2$.
If $|p_k|\ge 7$, then we have $-3m-n=\pm 1$ and $n=0$.
Hence we have no root for $(m, n)$.
If $|p_k|=5$, then we have $-3m-2n=\pm 1$.
Hence we have $(m, n)=(-1, 1)$ or $(1, -1)$.

\medskip

\noindent
(c) The case $c=|p_k|-2$.

\medskip

The lefthand side of (\ref{eq:lift1}) is 
$$\pm u^l(u^{|p_k|-1}-u^{|p_k|-2}-u+1)
=\mp u^{l+1}(2u^{|p_k|-3}+u^{|p_k|-4}+\cdots +u+2).$$
By (\ref{eq:deform}), we have $-m-3n=\pm 1$.
If $|p_k|\ge 7$, then we have $-n=\pm 2$ and $-m-n=\pm 1$.
Hence we have no root for $(m, n)$.
If $|p_k|=5$, then we have $-3m-2n=\pm 2$.
Hence we have no root for $(m, n)$.

\medskip

Therefore we have the result.

\bigskip

\noindent
(2) Suppose that $f_k=0$.

\medskip

The equation (\ref{eq:lenseq}) is simplified as
$$
\left(
\prod_{\scriptsize{j\ne k}
\atop
\scriptsize{1\le j\le \lambda}}
p_j
\right)
(\zeta_d^a-1)(\zeta_d^b-1)\doteq (\zeta_d-1)(\zeta_d^c-1).
$$
By taking the $d$-norm of both sides,
we have the result.
\qed

\bigskip

By (\ref{eq:Alex}), 
the Alexander polynomial of $M_3$ satisfies $f=f_k=1$ corresponding to 
Lemma \ref{lm:key2} (1), and
that of $M_{\lambda}$ ($\lambda \ge 4$) satisfies $f=f_k=0$ corresponding to 
Lemma \ref{lm:key2} (2).

\begin{figure}[htbp]
\begin{center}
\includegraphics[scale=0.8]{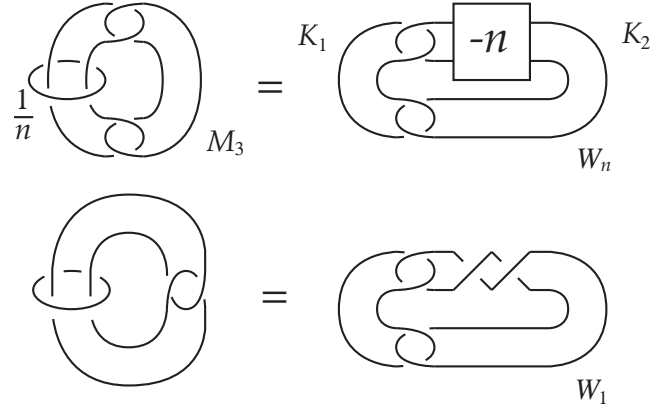}
\label{Wn}
\caption{$n$-twisted Whitehead link $W_n$}
\end{center}
\end{figure}

\begin{lm}\label{lm:KMS}{\rm \cite[Theorem 1.2]{KMS}}
Let $W_n\ (n\ne 0)$ be the $n$-twisted Whitehead link as in Figure 2
where the rectangle implies $(-n)$-full twists,
and $Y=(W_n; p_1/q_1, p_2/q_2)$
the result of a finite slope surgery along $W_n$.

\begin{enumerate}
\item[(a)]
If $|n|\ge 2$, then $Y$ is not a lens space.

\item[(b)]{\rm \cite{MP}}
Suppose that $n=\varepsilon=1$ or $-1$.
Then $Y$ is a lens space if and only if
one of the following holds:

\begin{enumerate}
\item[(1)]
$p_1/q_1=\varepsilon$ and $|\varepsilon p_2-6q_2|=1,$ or

\item[(2)]
$p_1/q_1=2\varepsilon$ and $|\varepsilon p_2-4q_2|=1,$ or

\item[(3)]
$p_1/q_1=3\varepsilon$ and $|\varepsilon p_2-3q_2|=1,$

\end{enumerate}
and the cases that indices of $(p_i, q_i)\ (i=1, 2, 3)$  are permuted.
Moreover, if (1), then $Y=L(p_2, 4\varepsilon q_2)$, 
if (2), then $Y=L(2p_2, \varepsilon (8q_2-p_2))$, and 
if (3), then $Y=L(3p_2, \varepsilon (3q_2-2p_2))$.

\end{enumerate}
\end{lm}

\section{Proof of Theorem \ref{th:MT1}}~\label{sec:3-comp}
Since the components of $M_3$ are interchangeable one another,
we may suppose $1\le |p_1|\le |p_2|\le |p_3|$.
By Lemma \ref{lm:key2} (1), we have $|p_1|=1$.

\medskip

\noindent
(i) The case $|p_3|\ge 5$.

\medskip

We set $p_1/q_1=\varepsilon$.
Since 
$(M_3; \varepsilon, p_2/q_2, p_3/q_3)
=(W_{\varepsilon}; p_2/q_2, p_3/q_3)$
and Lemma \ref{lm:KMS} (b), we have the result.

\bigskip

\noindent
(ii) The case $|p_3|\le 4$.

\medskip

We set $p_1/q_1=1/n$\ $(n\ne 0)$.
Since 
$(M_3; 1/n, p_2/q_2, p_3/q_3)
=(W_n; p_2/q_2, p_3/q_3)$
and Lemma \ref{lm:KMS} (a), we have $n=\pm 1$.
By the similar argument as (i), we have the result.
\qed

\section{Proofs of Theorem \ref{th:MT2}\ and Corollary \ref{co:propPR}}~\label{sec:many-comp}
We suppose that $L$ and $\overline{L}$ are the same links as in Section \ref{sec:key},
$\lambda \ge 4$, and $|p_i|\ge 2$ for fixed $i$.
By Lemma \ref{lm:key2} $(2),$ we may assume $p_j=1\ (j\ne i)$.
Let $\widehat{K}_i$ be the resulting knot by $1/q_j$-surgery along 
$K_j$\ $(j\ne i)$ of $L$, and
$\overline{L}_i$ the resulting $2$-component link by $1/q_j$-surgery along 
$K_j$\ $(j\ne i)$ of $\overline{L}$.
We set $Y_i=E_{\widehat{K}_i}$ and $\overline{Y}_i=E_{\overline{L}_i}$.
We note that $\widehat{K}_i$ and $\overline{L}_i$ are a knot and a link in $S^3$.
Then we have

\begin{eqnarray*}
H_1(Y_i) & = & 
\langle [m_1], [l_1]\ldots, [m_{\lambda}], [l_{\lambda}]\ |\ 
[m_j][l_j]^{q_j}=1, [l_j]=1\ (j\ne i), [l_i]=1\rangle \medskip\\
 & \cong & \langle [m_1], \ldots, [m_{\lambda}]\ |\ [m_j]=1\ (j\ne i)\rangle
\medskip\\
 & \cong & \langle [m_i]\ |\ -\rangle \cong \mathbb{Z},
\end{eqnarray*}

\begin{equation}\label{eq:barYi}
\begin{array}{ccl}
H_1(\overline{Y}_i) & = & 
\left\langle 
\begin{matrix}
[m_1], [l_1], \ldots, [m_{\lambda}], [l_{\lambda}], \hfill\\
[m_1'], [l_1']\ldots, [m_{\lambda}'], [l_{\lambda}'],\hfill\\
[m], [l] \hfill
\end{matrix}
\ \left|\
\begin{matrix}
[m_j']=[m_j][l_j]^{q_j}=1, \hfill\\
[l_j']=[l_j]^{-1}\ (j\ne i), \hfill\\
[l_j]=[m]\ (j=1, \ldots, \lambda),\hfill\\
[l]=[m_1]\cdots [m_{\lambda}] \hfill
\end{matrix}
\right.
\right\rangle
\medskip\\
& \cong & 
\langle 
[m_1], \ldots, [m_{\lambda}], [m]\ |\ 
[m_j][m]^{q_j}=1\ (j\ne i)
\rangle
\medskip\\
& \cong & \langle 
[m_i], [m]\ |\ -\rangle \cong \mathbb{Z}^2,
\end{array}
\end{equation}
and
\begin{equation}\label{eq:taubarYi}
\begin{matrix}
\tau(\overline{Y}_i)
& \doteq & {\displaystyle
{\mit \Delta}_{\overline{L}}
(t^{-q_1}, \ldots, t_i, \ldots, t^{-q_{\lambda}}, t)
\prod_{\scriptsize{j\ne i}
\atop
\scriptsize{1\le j\le \lambda}}
([l_j']-1)^{-1}}\hfill
\medskip\\
& \doteq &{\displaystyle
(t-1)g(t^{-q_1}, \ldots, t_i, \ldots, t^{-q_{\lambda}}, t)
\prod_{\scriptsize{j\ne i}
\atop
\scriptsize{1\le j\le \lambda}}
(t-1)^{-1}
}\hfill
\medskip\\
& = & {\displaystyle 
t_it^{-q_1-\cdots -\hat{q}_i-\cdots -q_{\lambda}}
}\hfill\medskip\\
&  & {\displaystyle 
+(t-1)(t_i-1)\sum_{
{\scriptstyle i\in I\subset I_{\lambda}}
\atop
{\scriptstyle 2\le s=|I| \le \lambda-1}}
\prod_{j\in I\setminus \{ i \}}\frac{t^{-q_j}-1}{t-1}h_I}\hfill\medskip\\
&  & {\displaystyle 
+(t-1)^2\sum_{
{\scriptstyle i\not \in I\subset I_{\lambda}}
\atop
{\scriptstyle 2\le s=|I| \le \lambda-1}}
\prod_{j\in I}\frac{t^{-q_j}-1}{t-1}h_I}\hfill\medskip\\
& & {\displaystyle +(t-1)(t_i-1)
\prod_{\scriptsize{j\ne i}
\atop
\scriptsize{1\le j\le \lambda}}
\frac{t^{-q_j}-1}{t-1}h}\hfill
\end{matrix}
\end{equation}
where $h_I$ and $h$ are the same as in Lemma \ref{lm:form}.

\medskip

\begin{lm}\label{lm:trivAlex}
${\mit \Delta}_{\widehat{K}_i}(t_i)\doteq 1$.
\end{lm}

\noindent
{\it Proof.}\ 
Let $\pi : \mathbb{Z}[t_i^{\pm 1}, t^{\pm 1}]\to \mathbb{Z}[t_i, t_i^{-1}]$
be a ring homomorphism defied by $\pi(t_i)=t_i$ and $\pi(t)=1$.
Then by (\ref{eq:barYi}) and (\ref{eq:taubarYi}), we have
$$
\tau(Y_i)\doteq 
\tau^{\pi}(\overline{Y}_i)(\pi([l])-1)^{-1}
=\frac{t_i}{t_i-1}.$$
By Proposition \ref{pr:torEx} (3), we have the result.
\qed

\begin{theo}\label{th:Ni}{\rm \cite[Corollary 1.3]{Ni}}
If a knot $K$ in $S^3$ yields a lens space, then $K$ is fibered.
\end{theo}

By this theorem, 
existence of lens surgery along $\widehat{K}_i$
is equivalent to triviality of $\widehat{K}_i$.

\begin{theo}\label{th:Gabai}{\rm \cite{Gb1}}
Let $F$ be a Seifert surface of a link $L$
which is a Murasugi sum of two surfaces $F_1$ and $F_2$.
Then $F$ is a minimal genus Seifert surface of $L$ if and only if 
both $F_1$ and $F_2$ are minimal genus Seifert surfaces of links.
\end{theo}

\noindent
{\it Proof of Theorem \ref{th:MT2}.}\ 

\medskip

\noindent
(1) The case $\lambda=4$.

\medskip

By the symmetry of $M_4$ as in Section \ref{sec:intro},
it is sufficient to prove that $\widehat{K}_1$ is non-trivial.
Let $F$ be a genus $1$ Seifert surface of $M_4$ as in Figure 3.
Then $F$ is a plumbing of $F_1$ and $F_2$ in Figure 3.
In Figure 3, the number in a rectangle implies the number of full twists.
Since $\partial F_1$ is a non-trivial $(2, -2q_2)$-torus link, 
$F_1$ is a minimal genus Seifert surface.
Since $\partial F_1=K_1'\cup K_2'$ is a parallel of 
a non-trivial $2$-bridge knot $C(-2q_3, 2q_4)$,
$F_2$ is a minimal genus Seifert surface.
By Theorem \ref{th:Gabai},
$F$ is a minimal genus Seifert surface, and hence
$\widehat{K}_1$ is a non-trivial knot.

\bigskip

\noindent
(2) The case $\lambda \ge 5$.

\medskip

Let $Y$ be the result of Dehn surgery along $M_{\lambda}$
whose surgery coefficient of $K_j$ ($j\ne i$) is $1/q_j$ and
that of $K_i$ is $p/q$.
We show that $Y$ includes an essential torus.
Let $T$ be a standard torus embedded in $S^3$ which divides $S^3$
into two solid tori $\Sigma_1$ and $\Sigma_2$.
We can take $T$ such that
$\Sigma_1$ includes $L_1=K_1\cup K_2$, and
$\Sigma_2$ includes $L_2=K_3\cup \cdots \cup K_{\lambda}$ (Figure 4).
We show that both the result of Dehn surgery along $L_1$ in $\Sigma_1$,
which is denoted by $W_1$,
and the result of Dehn surgery along $L_2$ in $\Sigma_2$,
which is denoted by $W_2$, are not solid tori.

\medskip

Suppose that $W_1$ is a solid torus.
Then in Theorem \ref{th:MT1},
there exists one of three sequences such that a surgery coefficients
of some component can be any rational number.
Hence $W_1$ is not a solid torus.

\medskip

Suppose that $W_2$ is a solid torus.
Since the complement of $L_2$ in $\Sigma_2$
is homeomorphic to $M_{\lambda-1}$ with preserving
the peripheral structures of the components except $\partial \Sigma_2$,
it contradicts the assumption of induction.
Hence $W_2$ is not a solid torus.
\qed

\begin{figure}[htbp]
\begin{center}
\includegraphics[scale=0.8]{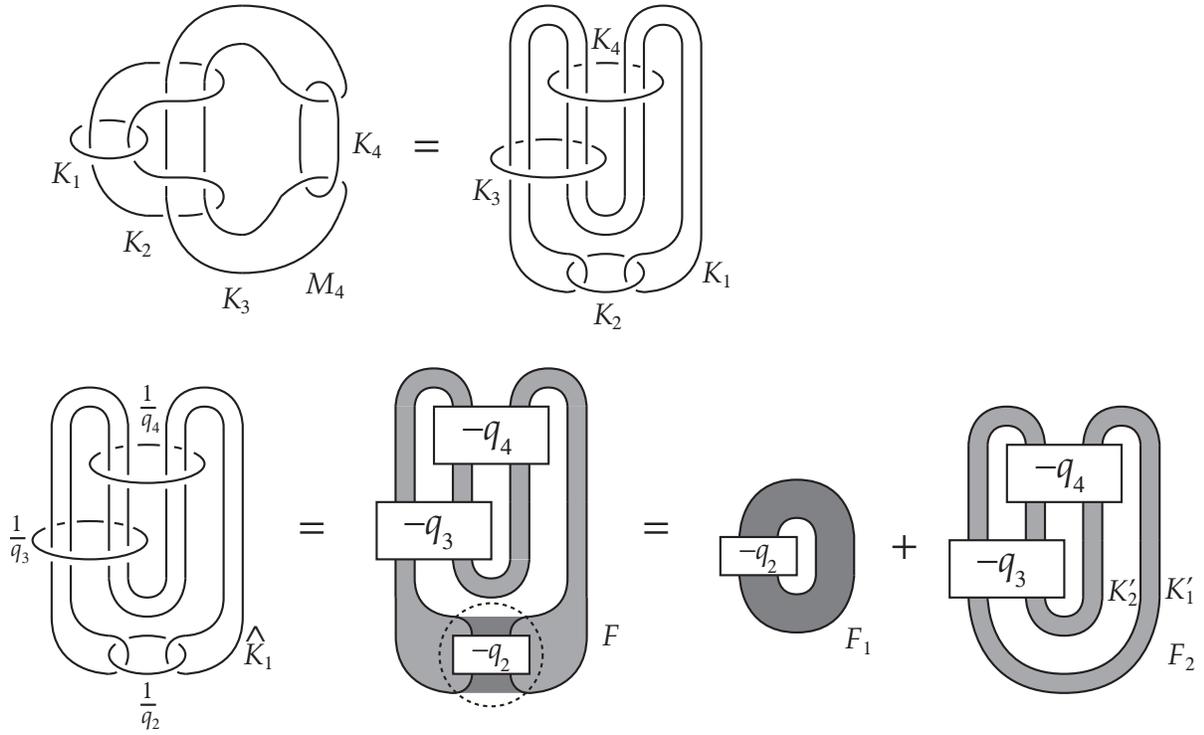}
\label{plumb}
\caption{minimal genus Seifert surface for 
$\widehat{K}_1$}
\end{center}
\end{figure}

\begin{figure}[htbp]
\begin{center}
\includegraphics[scale=0.75]{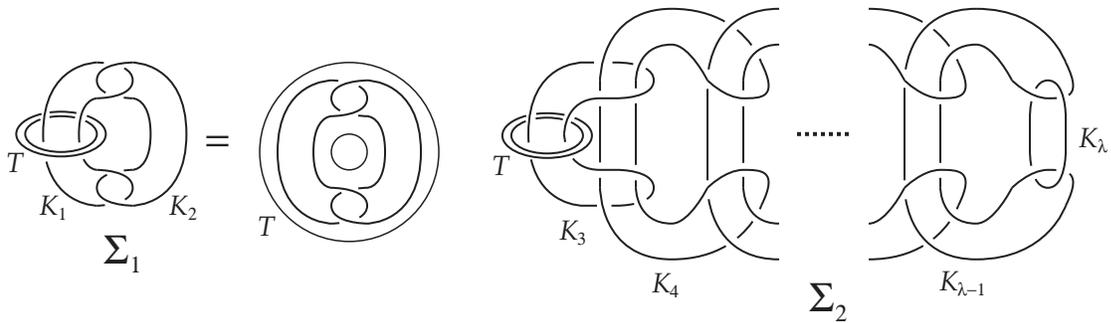}
\label{essential}
\caption{essential torus $T$, and solid tori $\Sigma_1$ and $\Sigma_2$}
\end{center}
\end{figure}

\begin{re}\label{rm:CW}
{\rm 
\begin{enumerate}
\item[(1)]
In the proof of Theorem \ref{th:MT2},
the case $\lambda=4$ can also be proved in a similar way as
the case $\lambda \ge 5$.
It will be discussed again in Subsection \ref{ssec:geogen}.

\item[(2)]
We cannot prove Theorem \ref{th:MT2} 
by only using the Reidemeister torsions and
the Casson-Walker invariants \cite{Les, Wa}.
Hence an invariant induced from the knot Floer homology
due to L.~Nicolaescu \cite{Nc} also cannot prove it.
However the knot Floer homology itself can prove it by Ni's theorem above.

\end{enumerate}}
\end{re}

\noindent
{\it Proof of Corollary \ref{co:propPR}.}\ 
Suppose that $H_1(Y)\cong \{0\}$ (resp.\ $H_1(Y)\cong \mathbb{Z}$).
There exists $i\in \{1, 2, \ldots, \lambda\}$
such that a surgery coefficient of $K_j$ with $j\ne i$ is $1/q_j$
and that of $K_i$ is $1/q_i$ (resp.\ $0$).
Then since $\widehat{K}_i$ is a $2$-bridge knot of genus $1$ for $\lambda=3$
and $\widehat{K}_i$ is non-trivial by the proof of Theorem \ref{th:MT2} for $\lambda \ge 4$,
$Y$ is not $S^3$ (resp.\ not $S^1\times S^2$) 
by the affirmative answer of Property P conjecture (resp.\ Property R conjecture)
due to P.~Kronheimer and T.~Mrowka \cite{KM}
(resp.\ D.~Gabai \cite{Gb2}).
\qed

\section{Generalization to Brunnian type links
and toroidal Brunnian links}~\label{sec:gen}
We generalize Theorem \ref{th:MT1}, Theorem \ref{th:MT2} and Corollary \ref{co:propPR} 
for the cases of {\it Brunnian type links} and 
toroidal Brunnian type links
which are Brunnian links including essential tori in the link complement.
Moreover we characterized 
toroidal Brunnian links and toroidal Brunnian type links.

\subsection{Generalization to Brunnian type links}~\label{ssec:alggen}
A $\lambda$-component link $L=K_1\cup \cdots \cup K_{\lambda}$ in $S^3$ 
with $\lambda \ge 2$ is a {\it Brunnian type link} if $L$ is an algebraically split link
such that every component is an unknot
and the Alexander polynomial of every proper sublink with at least two components is zero.
Let $\overline{L}=L\cup K$ be the same link as in Subsection \ref{ssec:asplit}.
Then we have the Alexander polynomial of $\widehat{K}_i$ 
in Section \ref{sec:many-comp} by a similar way as
(\ref{eq:taubarYi}) and the proof of Lemma \ref{lm:trivAlex} : 
\begin{equation}\label{eq:hatKi}
{\displaystyle
{\mit \Delta}_{\widehat{K}_i}(t_i)\doteq
t_i+(-1)^{\lambda-1}\left(\prod_{\scriptsize{j\ne i}
\atop
\scriptsize{1\le j\le \lambda}}q_j
\right)f_i(t_i)(t_i-1)^2}.
\end{equation}
We note that by Remark \ref{re:g}, the sign of $f_i(t_i)$ is uniquely determined.
By  (\ref{eq:hatKi}) 
and the affirmative answer of Property P conjecture (resp.\ Property R conjecture)
due to P.~Kronheimer and T.~Mrowka \cite{KM}
(resp.\ D.~Gabai \cite{Gb2}), we have the following.

\begin{theo}\label{th:MT3}
We suppose the situation above
(i.e.\ $L$ is a $\lambda$-component Brunnian type link).
\begin{enumerate}
\item[(1)]
If the result of a finite slope surgery along $L$ is $S^3$, then
$f_i=0$ for all $i\in \{1, 2, \ldots, \lambda\}$.

\item[(2)]
If the result of a finite slope surgery along $L$ is $S^1\times S^2$, then
$f_i=0$ for some $i\in \{1, 2, \ldots, \lambda\}$.

\end{enumerate}
\end{theo}

Let $L=K_1\cup \cdots \cup K_{\lambda}$ be a $\lambda$-component link in $S^3$
and $W_1=K_1'\cup K_2'$ the Whitehead link (cf.\ Figure 2).
Since $K_2'$ is an unknot, the complement of $K_2'$ is a solid torus $\Sigma$ and
we may regard $W_1$ as a knot $K_1'$ in $\Sigma$.
If we attach $\Sigma$ to $\partial N(K_i)$ by identifying 
the preferred longitude of $K_2'$ with the meridian of $K_i$
and the meridian of $K_2'$ with the preferred longitude of $K_i$,
then we obtain a new $\lambda$-component link $L'$ in $S^3$ such that
$K_i$ is replaced with $K_1'$.
We call $L'$ a {\it Whitehead double} of $L$ on $K_i$.

\begin{ex}\label{ex:Brunnian}
{\rm
Let $L=K_1\cup K_2\cup K_3$ be a $3$-component link 
obtained by taking a Whitehead double on one component of the Whitehead link
and taking $2$-parallel of the other component of the Whitehead link.
We suppose that $K_1$ is the Whitehead double part and
$K_2\cup K_3$ is the $2$-parallel part.
Then we can easily see that $L$ is a Brunnian type link,
but is not a Brunnian link satisfying ${\mit \Delta}_L(t_1, t_2, t_3)=0$.
Since $(L; \emptyset, 1/q, -1/q)$ is an unknot for every integer $q$,
$L$ has not Property FP$'$ and Property $FR$.
}
\end{ex}

\begin{figure}[htbp]
\begin{center}
\includegraphics[scale=0.8]{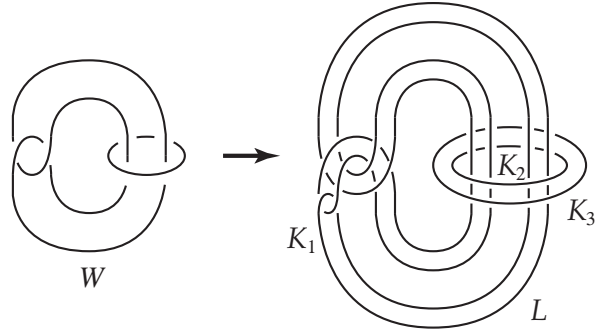}
\label{ex}
\caption{Brunnian type link $L=K_1\cup K_2\cup K_3$}
\end{center}
\end{figure}

We do not have any counterexamples for the following problems.

\begin{qu}\label{qu:Brunnian}
\begin{enumerate}
\item[(1)]
Does any non-split $\lambda$-component Brunnian type link $L$
without Property FP$'$ or Property FR satisfy
${\mit \Delta}_L(t_1, \ldots, t_{\lambda})=0$ ?

\item[(2)]
Does any non-split anannular Brunnian type link have Property FP$'$ or Property $FR$ ?

\end{enumerate}
\end{qu}

Corollary \ref{co:propPR} is a partial answer for Question \ref{qu:Brunnian} (2).

\medskip

About lens surgeries along a knot, 
P.~Ozsv\'ath and Z.~Szab\'o \cite{OS} obtained the following
by the knot Floer homology.

\begin{theo}\label{th:OS}
{\rm (Ozsv\'ath-Szab\'o \cite[Corollary 1.3]{OS})}
Let $K$ be a knot in $S^3$.
If $K$ yields a lens space, then the Alexander polynomial of $K$ is of the following form : 
$${\mit \Delta}_K(t)\doteq
(-1)^m+\sum_{k=1}^m(-1)^{k-1}(t^{n_k}+t^{-n_k})\quad
(n_1>n_2>\cdots >n_m\ge 1)$$
where $n_k\ (k=1, \ldots, m)$ is an integer.
\end{theo}

The author and Y.~Yamada \cite{KY} obtained the following
by the Reidemeister torsion.

\begin{theo}\label{th:KY}
{\rm \cite{KY}}
Let $K$ be a knot in an integral homology $3$-sphere.
If $(K; p/q)\ (p\ge 2)$ is a lens space, then 
the Alexander polynomial of $K$ is of the following form : 
$${\mit \Delta}_K(t)\doteq
\frac{(t^{rs}-1)(t-1)}{(t^r-1)(t^s-1)}\quad \left(\mathrm{mod}\ \! \frac{t^p-1}{t-1}\right)$$
where $\gcd(r, s)=1$ and $qrs\equiv \pm 1\ (\mathrm{mod}\ \! p)$.
\end{theo}

By combining Theorem \ref{th:OS} and Theorem \ref{th:KY}, 
M.~Tange \cite{Ta} obtained the following.
M.~Hedden and L.~Watson \cite{HW} obtained the same result independently.

\begin{theo}\label{th:Tange}
{\rm (Tange \cite[Theorem 1.3]{Ta}, Hedden and Watson \cite[Corollary 9]{HW})}
Under the same situation as Theorem \ref{th:OS}, we have $n_1-n_2=1$.
That is, the Alexander polynomial of $K$ is of the following form : 
$${\mit \Delta}_K(t)\doteq
t^{n_1}-t^{n_1-1}+\cdots -t^{-n_1+1}+t^{-n_1}.$$
\end{theo}

Theorem \ref{th:Tange} implies that
the trace, which is the sum of the roots, of ${\mit \Delta}_K(t)$ is $1$.

\begin{theo}\label{th:MT4}
We suppose the same situation in Theorem \ref{th:MT3}.
If the result of a finite slope surgery along $L$ is a lens space, then
we have the following:
\begin{enumerate}
\item[(1)]
For some $i\in \{1, 2, \ldots, \lambda\}$ and 
every $j\in \{1, 2, \ldots, \lambda\}\setminus \{i\}$,
we have $|q_j|=1$.

\item[(2)]
For the same $i\in \{1, 2, \ldots, \lambda\}$ in (1),
${\mit \Delta}_{\widehat{K}_i}(t_i)$ satisfies
the conditions in Theorem \ref{th:OS}, Theorem \ref{th:KY} and Theorem \ref{th:Tange}.

\end{enumerate}
\end{theo}

By (\ref{eq:hatKi}), Theorem \ref{th:MT3} and Theorem \ref{th:MT4}, we have :

\begin{co}\label{co:const}
Under the same situation in Theorem \ref{th:MT4},
if $f_i(t_i)$ is constant, then we have
(i) $f_i=0$ and $\widehat{K}_i$ is an unknot, or 
(ii) 
$$f_i(t_i)=(-1)^{\lambda-1}\left(\prod_{\scriptsize{j\ne i}
\atop
\scriptsize{1\le j\le \lambda}}q_j
\right) =1\quad 
(\mbox{i.e.\ $|q_j|=1$ for every $j\in \{1, 2, \ldots, \lambda\}\setminus \{i\}$}),$$
$${\mit \Delta}_{\widehat{K}_i}(t_i)\doteq t_i^2-t_i+1,$$
and $\widehat{K}_i$ is the trefoil.
\end{co}

\subsection{Generalization to toroidal Brunnian type links}~\label{ssec:geogen}
A link $L$ is {\it toroidal} if there exists an essential torus in the complement of $L$.
A link $L$ is {\it atoroidal} if $L$ is not toroidal.

\medskip

Let $L=K_1\cup \cdots \cup K_p$ be a $p$-component link in $S^3$, 
and $L'=K_1'\cup \cdots \cup K_q'$ a $q$-component link in $S^3$.
Suppose that for some $j$, the $j$-th component $K_j'$ of $L'$ is an unknot.
Then the complement of $K_j'$ is a solid torus, and
we denote it by $\Sigma'$.
If we attach $\Sigma'$ to the complement of $L$
by identifying the preferred longitude of $K_j'$ with the meridian of $K_i$ 
and the meridian of $K_j'$ with the preferred longitude of $K_i$, 
then we obtain a new $(p+q-2)$-component link in $S^3$.
We call the link the {\it satellite} of $L$ along $K_i$ 
with {\it pattern} $(L', K_j')$ or $L'$ of type $(p, q)$.
We denote it by $S_{ij}(L, L')$ or $S(L, L')$ simply.
If we set the attaching map $f : \Sigma' \to S^3$, then
$S_{ij}(L, L')=(L\setminus K_i)\cup f(L'\setminus K_j')$.
If $L'$ is the Whitehead link, then $S_{ij}(L, L')$ is the {\it Whitehead double} of $L$, 
if $L'$ is the Borromean rings, then $S_{ij}(L, L')$ is the {\it Bing double} of $L$, 
and if $L$ is the Hopf link, then $S_{ij}(L, L')=L'$.
In particular, $M_3$ is the Bing double of the Hopf link, and
$M_{\lambda}\ (\lambda \ge 4)$ is the Bing double of $M_{\lambda-1}$.

\medskip

A $\lambda$-component link $L=K_1\cup \cdots \cup K_{\lambda}$ in $S^3$
is a {\it semi-Brunnian link} if it is a Brunnian link for $\lambda \ge 3$,
every component of it is an unknot for $\lambda=2$, and
it is an arbitrary knot for $\lambda=1$.
Then we have the following.

\begin{lm}\label{lm:construct}
\begin{enumerate}
\item[(1)]
Let $L=K_1\cup \cdots \cup K_p$ be a $p$-component semi-Brunnian link,
and $L'=K_1'\cup \cdots \cup K_q'$ a $q$-component Brunnian link 
such that $p=1$ and $q \ge 3$, or $p \ge 2$ and $q \ge 2$.
Then $S_{ij}(L, L')$ is a $(p+q-2)$-component Brunnian link in $S^3$
for any $i\in \{1, 2, \ldots, p\}$ and $j\in \{1, 2, \ldots, q\}$.

\item[(2)]
In (1), $S_{ij}(L, L')$ is non-split if and only if both $L$ and $L'$ are non-split.

\item[(3)]
In (1), if $L$ is non-split atoroidal and 
is neither the unknot nor the Hopf link, 
and $L'$ is non-split atoroidal, then
$S_{ij}(L, L')$ is a $(p+q-2)$-component non-split toroidal Brunnian link in $S^3$
for any $i\in \{1, 2, \ldots, p\}$ and $j\in \{1, 2, \ldots, q\}$ with a unique essential torus.

\end{enumerate}
\end{lm}

\noindent
{\it Proof.}\ 
(1) Let $T$ be an embedded solid torus in $S^3$ such that
$T$ separates $S_{ij}(L, L')$ into two parts
$L\setminus K_i$ and $L'\setminus K_j'$, and
$T$ separates $S^3$ into two parts $\Sigma$ and $\Sigma'$
where $\Sigma$ is the complement of $K_i$ and
$\Sigma'$ is the complement of $K_j'$.
We note that $\Sigma'$ is a solid torus, and $\Sigma$ is also a solid torus if $p\ge 2$
and is never a solid torus if $p=1$.

\medskip

Suppose that $p\ge 2$.
For $r\in \{1, 2, \ldots, p\}\setminus \{ i\}$, 
we set $L_r=L\setminus K_r$.
If $p=2$, then $L_r$ is the empty set.
Then we have 
$S_{ij}(L, L')\setminus K_r=S_{ij}(L_r, L')$.
Since $L_r$ is a $(p-1)$-component trivial link
and $L'\setminus K_j'$ is a $(q-1)$-component trivial link, 
$S_{ij}(L, L')\setminus K_r$ is a $(p+q-3)$-component trivial link.

\medskip

For $s\in \{1, 2, \ldots, q\}\setminus \{ j\}$, 
we set $L_s'=L'\setminus K_s'$.
Then we have 
$S_{ij}(L, L')\setminus K_s'=S_{ij}(L, L_s')$.
Since $L_s'$ is a $(q-1)$-component trivial link
and $L\setminus K_i$ is a $(p-1)$-component trivial link,
$S_{ij}(L, L')\setminus K_s'$ is a $(p+q-3)$-component trivial link.

\bigskip

\noindent
(2) We use the same notations as in the proof of (1).
If $L$ or $L'$ is split, then it is easy to see that $S_{ij}(L, L')$ is also split.

\medskip

Suppose that $S_{ij}(L, L')$ is split.
Let $F$ be an embedded $2$-sphere in the complement of $S_{ij}(L, L')$
(i.e.\ $F\cap S_{ij}(L, L')=\emptyset$)
which separates $S_{ij}(L, L')$ into two non-empty links.
If $T\cap F=\emptyset$, then $L$ or $L'$ is split.
We show that the case $T\cap F\ne \emptyset$ does not occur.
We suppose $T\cap F\ne \emptyset$, and $T$ intersects with $F$ transverselly.
Then $T\cap F$ consists of disjoint simple loops.
We take an innermost loop $\gamma$ of $T\cap F$ on $F$
and the disk $B$ bounded by $\gamma$ on $F$.
Then $B$ is properly embedded in $\Sigma$ or $\Sigma'$.
We only show the case that $B$ is properly embedded in $\Sigma$.
If $\gamma$ bounds a disk $B'$ on $T$, then
$B\cup B'$ bounds a $3$-ball $W$ in $\Sigma$.
Since $L$ is non-split, $W\cap L=\emptyset$ 
and $\gamma$ can be removed by local isotopy of $F$ along $W$.
If $\gamma$ is not null-homotopic on $T$, then
$\gamma$ is a preferred longitude of $K_i$.
Since $L$ is non-split, 
$B\cap L=B\cap (L\setminus K_i)\ne \emptyset$.
It contradics the assumption that $F\cap S_{ij}(L, L')=\emptyset$.
Hence the case $T\cap F\ne \emptyset$ does not occur.

\bigskip

\noindent
(3) We use the same notations as in the proof of (1).
By (1) and (2), 
$S_{ij}(L, L')$ is a $(p+q-2)$-component non-split Brunnian link in $S^3$
for any $i\in \{1, 2, \ldots, p\}$ and $j\in \{1, 2, \ldots, q\}$.

\medskip

Let $E$ be the complement of $S_{ij}(L, L')$.
Since $L$ is neither the unknot nor the Hopf link,
$T$ neither bound a solid torus nor is boundary parallel in $E$.
Therefore $T$ is an essential torus in $E$.

\medskip

We show that $T$ is a unique essential torus in $E$.
Suppose that there exists another essential torus $T'$ in $E$
(i.e.\ $T'\cap S_{ij}(L, L')=\emptyset$ and
$T'$ is not ambient isotopic to $T$ in $E$).
If $T\cap T'=\emptyset$, then $T'$ is an essential torus
in the complement of $L$ or $L'$.
By the assumption, the case does not occur.

\medskip

We may suppose that $T\cap T'\ne \emptyset$, and
$T'$ intersects with $T$ transverselly.
Then $T\cap T'$ consists of disjoint simple loops, and
$T'$ is divided into $T'\cap \Sigma$ and $T'\cap \Sigma'$ by $T$.
We take a loop $\gamma$ of $T\cap T'$.
If $\gamma$ bounds a disk $B'$ on $T'$, then
$\gamma$ also bounds a disk $B$ on $T$, 
or $\gamma$ is a preferred longitude of $K_i$ or $K_j'$.
If $\gamma$ bounds a disk $B$ on $T$, then
$B\cup B'$ bounds a $3$-ball $W$ in $\Sigma$ or $\Sigma'$.
Since both $L$ and $L'$ are non-split, 
we have $W\cap L=W\cap L'=\emptyset$, and
$\gamma$ can be removed by local isotopy of $T'$ along $W$.
If $\gamma$ is a preferred longitude of $K_i$, then
$B'$ is properly embedded in $\Sigma$.
Since $L$ is non-split, we have $B'\cap L=B'\cap (L\setminus K_i)\ne \emptyset$.
It contradicts the assumption $T'\cap S_{ij}(L, L')=\emptyset$.
Hence this case does not occur.
We can show that $\gamma$ is not a preferred longitude of $K_j'$
in a similar way.
Therefore we may assume that $T\cap T'$ consists of
parallel essential simple loops on $T'$, and
the connected components of both $T'\cap \Sigma$ and $T'\cap \Sigma'$
are annuli.
The set $T\cap T'$ is also the set of parallel essential simple loops on $T$.
We take a connected component $A'$ of $T'\cap \Sigma$ (resp. $T'\cap \Sigma'$),
which is an annulus.
Then $\partial A'$ bounds an annulus $A$ on $T$ (resp. $T'$), and
$A\cup A'$ is a torus.
By a slight isotopy, the torus $A\cup A'$ can be regarded as 
an embedded torus in $\mathrm{Int}(\Sigma)$ (resp. $\mathrm{Int}(\Sigma')$).
Since $\Sigma$ (resp. $\Sigma'$) is atoroidal, 
the torus bounds a solid torus which is included by a regular neighborhood of
$\partial \Sigma$ (resp. $\partial \Sigma'$).
We can remove $\partial A'$ by local isotopy of $T'$ along the solid torus.
Therefore we can isotope $T'$ such that $T\cap T'=\emptyset$,
and the case does not occur ($T'$ cannot exist).
\qed

\bigskip

We show the converse of Lemma \ref{lm:construct} (3)
for the case $p \ge 2$ and $q \ge 2$.

\begin{theo}\label{th:construct}
For a $\lambda$-component non-split toroidal Brunnian link 
$\widehat{L}$ with $\lambda \ge 2$, 
if there exists an essential torus $T$ in the complement of $\widehat{L}$
such that $T$ decomposes $\widehat{L}$ into 
a $(p-1)$-component sublink and a $(q-1)$-component sublink
with $p \ge 2$ and $q \ge 2$, then there exist 
a $p$-component non-split semi-Brunnian link $L$
and a $q$-component non-split Brunnian link $L'$
such that $\widehat{L}=S(L, L')$ with $\lambda=p+q-2$.
\end{theo}

\noindent
{\it Proof.}\ 
Let $E$ be the complement of $\widehat{L}$, and $T$ an essential torus in $E$.
Then $T$ decomposes $S^3$ into two parts $\Sigma$ and $\Sigma'$, and
one of $\Sigma$ and $\Sigma'$ is a solid torus.
We suppose that $\Sigma'$ is a solid torus, and
that $\Sigma$ includes $p'$ components of $\widehat{L}$
and $\Sigma'$ includes $q'$ components of $\widehat{L}$.
Then $\lambda=p'+q'$.
We denote the core knot of $\Sigma'$ by $K$, 
a meridian of $\Sigma'$ on $\partial \Sigma'$ by $K'$, 
the $p'$-component link included in $\Sigma$ by $K_1\cup \cdots \cup K_{p'}$, and
the $q'$-component link included in $\Sigma'$ by $K_1'\cup \cdots \cup K_{q'}'$.
We take
$L=K_1\cup \cdots \cup K_{p'}\cup K$ and
$L'=K_1'\cup \cdots \cup K_{q'}'\cup K'$.
Then $p=p'+1$, $q=q'+1$, and $\widehat{L}=S_{pq}(L, L')$.
We show that $L$ is a $p$-component non-split semi-Brunnian link
and $L'$ is a $q$-component non-split Brunnian link.

\medskip

Both $\Sigma$ and $\Sigma'$ include at least one component of $\widehat{L}$.
Firstly we show that $L$ is a semi-Brunnian link.
Since $\widehat{L}$ is a Brunnian link, 
$L\setminus K=\widehat{L}\setminus (L'\setminus K')$ is a trivial link.
For $i\in \{1, \ldots, p'\}$, 
we suppose that $K$ is a non-split component or a non-trivial component
in $L\setminus K_i$.
Then a longitude of $K$ does not represent the trivial element
in the fundamental group of the complement of $L\setminus (K_i\cup K)$.
We note that the longitude of $K$ can be situated on
$T=\partial \Sigma$, and we denote it by $l$.
The meridian of $\Sigma'$, which is a longitude of $K'$, is also on $T=\partial \Sigma'$, 
and we denote it by $l'$.
For the complement of $\widehat{L}\setminus K_i$, 
since both $l$ and $l'$ are essential in $\Sigma$ and $\Sigma'$ respectively,
$T$ is still essential.
It contradicts that $\widehat{L}\setminus K_i$ is a trivial link.
Hence $K$ is an unknot and a split component in $L\setminus K_i$,
and $L\setminus K_i$ is a trivial link.
Therefore $L$ is a semi-Brunnian link.

\medskip

Secondly we show that $L'$ is a Brunnian link.
Since $\widehat{L}$ is a Brunnian link, 
$L'\setminus K'=\widehat{L}\setminus (L\setminus K)$ is a trivial link.
For $j\in \{1, \ldots, q'\}$, 
we suppose that $L'\setminus K_j'$ is not a trivial link.
Then since $L'\setminus (K_j'\cup K')$ is a trivial link,
a longitude of $K'$ does not represent the trivial element
in the fundamental group of the complement of $L'\setminus (K_j'\cup K')$.
For the complement of $\widehat{L}\setminus K_j'$, 
since $l$ is essential in $\Sigma$ and $l'$ is essential in $\Sigma'$,
$T$ is still essential.
It contradicts that $\widehat{L}\setminus K_j'$ is a trivial link.
Hence $L'\setminus K_j'$ is a trivial link, and $L'$ is a Brunnian link.
\qed

\bigskip

For the case $p=1$ and $q \ge 3$,
the converse of Lemma \ref{lm:construct} (3) does not hold.

\begin{ex}\label{ex:rubber}
{\rm
Let $R_{\lambda}=K_1\cup \cdots \cup K_{\lambda}$
be a $\lambda$-component rubberband Brunnian link \cite{Ba} with $\lambda \ge 3$, 
and $\overline{R}_{\lambda}=R_{\lambda}\cup K'$ as in Figure 6,
where $K_{\lambda+1}=K_1$.
Then $R_{\lambda}$ is a non-split atoroidal Brunnian link, and
$\overline{R}_{\lambda}$ is not a Brunnian link but an algebraically split link.
For a non-trivial knot $K$,
$\widehat{L}=S_{1,\lambda+1}(K, \overline{R}_{\lambda})$ is a non-split toroidal Brunnian link.
\begin{figure}[htbp]
\begin{center}
\includegraphics[scale=0.8]{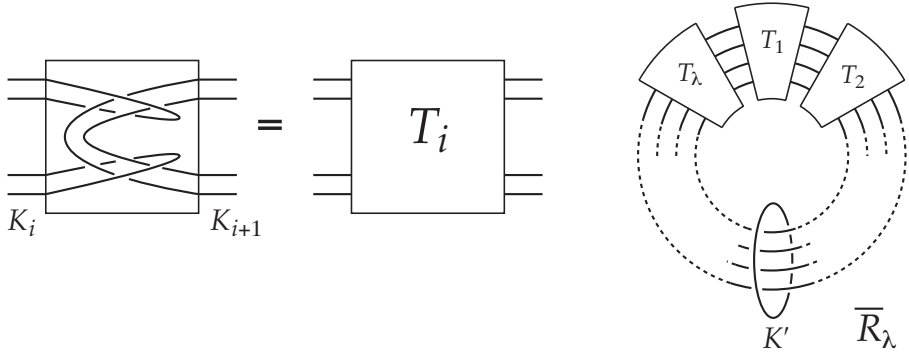}
\label{rubber}
\caption{a $\lambda$-component rubberband Brunnian link $R_{\lambda}$
and $\overline{R}_{\lambda}=R_{\lambda}\cup K'$}
\end{center}
\end{figure}
}
\end{ex}

\begin{re}\label{re:hyperbolic}
{\rm
Since a non-split Brunnian link is prime and anannular, 
a non-split atoroidal Brunnian link is hyperbolic.
There exists infinitely many hyperbolic Brunnian links \cite{Kn}.
For example, a rubberband Brunnian link $R_{\lambda}$
(see also Example \ref{ex:rubber} and Figure 6) is hyperbolic.}
\end{re}

A $\lambda$-component link $L=K_1\cup \cdots \cup K_{\lambda}$ in $S^3$
with $\lambda \ge 2$ is a {\it semi-Brunnian type link} if for $\lambda \ge 3$,
$L$ is algebraically split link such that for fixed $i\in \{1, 2, \ldots, \lambda\}$,
$K_r$ ($r\in \{1, 2, \ldots, \lambda\} \setminus \{i\}$) is an unknot 
and the Alexander polynomial of 
every proper sublink with at least two components is zero, and
for $\lambda =2$ and fixed $i\in \{1, 2\}$,
$K_r$ ($r\in \{1, 2\} \setminus \{i\}$) is an unknot.
Then $K_i$ is called the {\it characteristic component} of $L$.

\begin{theo}\label{th:Brunnian}
Let $L=K_1\cup \cdots \cup K_p$ be a $p$-component semi-Brunnian type link with $p \ge 2$,
$K_i$ the characteristic component of $L$,
and $L'=K_1'\cup \cdots \cup K_q'$ a $q$-component Brunnian type link with $q \ge 2$.
Then $\widehat{L}=S_{ij}(L, L')$ is a $(p+q-2)$-component Brunnian type link in $S^3$
for any $j\in \{1, 2, \ldots, q\}$.
\end{theo}

To prove Theorem \ref{th:Brunnian},
we prepare a lemma.

\begin{lm}\label{lm:BrunnianAlex}
In the same situation as in Theorem \ref{th:Brunnian},
we assume that $i=p$ and $j=q$.
Then we have the following:
\begin{enumerate}
\item[(1)]
$${\mit \Delta}_L(t_1, \ldots, t_p)\doteq
\left\{
\begin{array}{cl}
{\displaystyle
\frac{(t_1t_2)^k-1}{t_1t_2-1}{\mit \Delta}_{K_2}(t_2)
+(t_1-1)(t_2-1)f(t_1, t_2)} & (p=2),\medskip\\
(t_1-1)\cdots (t_p-1)f(t_1, \ldots, t_p) & (p\ge 3),
\end{array}
\right.$$
where $k=\mathrm{lk}(K_1, K_2)$ for $p=2$,
and
$${\mit \Delta}_{L'}(t_1', \ldots, t_q')\doteq
(t_1'-1)\cdots (t_q'-1)f'(t_1', \ldots, t_q'),$$
where $f(t_1, t_2)\in \mathbb{Z}[t_1^{\pm 1}, t_2^{\pm 1}]$, 
$f(t_1, \ldots, t_p)\in \mathbb{Z}[t_1^{\pm 1}, \ldots, t_p^{\pm 1}]$, and
$f'(t_1', \ldots, t_q')\in \mathbb{Z}[t_1'^{\pm 1}, \ldots, t_q'^{\pm 1}]$.

\item[(2)]
If $p=2$ and $k\ne 0$, then we have 
$${\mit \Delta}_{\widehat{L}}(t_1, t_1', \ldots, t_{q-1}')
\doteq
\frac{(t_1^k-1)^2}{t_1-1}(t_1'-1)\cdots (t_{q-1}'-1)f'(t_1', \ldots, t_{q-1}', t_1^k).$$

\item[(3)]
If $p=2$ and $k=0$, or $p\ge 3$, then we have 
$${\mit \Delta}_{\widehat{L}}(t_1, \ldots, t_{p-1}, t_1', \ldots, t_{q-1}')
\doteq 0.$$

\end{enumerate}

\end{lm}

\noindent
{\it Proof.}\ 
(1) By applying Lemma \ref{lm:Torres} repeatedly, we have the result.

\medskip

\noindent
(2) By (1) and Lemma \ref{lm:surgery1}, we have the result.

\medskip

\noindent
(3) We add one component $K$ to $L$ such that
$K_i\cup K$ ($i=1, \ldots, p$) is the connected sum of $K_i$ and the Hopf link, 
where $\mathrm{lk}(K_i, K)=1$, and 
we denote by $\overline{L}=L\cup K$
and $\widehat{\overline{L}}=\widehat{L}\cup K$.
Then by (1), (\ref{eq:Torres1}) and Lemma \ref{lm:surgery1}, we have
$$\begin{array}{ccl}
{\mit \Delta}_{\widehat{\overline{L}}}(t_1, \ldots, t_{p-1}, t_1', \ldots, t_{q-1}', t)
& \doteq &
{\mit \Delta}_{\overline{L}}(t_1, \ldots, t_{p-1}, 1, t)
{\mit \Delta}_{L'}(t_1', \ldots, t_{q-1}', t)\medskip\\
& \doteq &
(t-1)^2g(t_1, \ldots, t_{p-1}, 1, t)\medskip\\
& & \times (t_1'-1)\cdots (t_{q-1}'-1)f'(t_1', \ldots, t_{q-1}', t),
\end{array}
$$
where $g(t_1, \ldots, t_{p-1}, 1, t)\in \mathbb{Z}[t_1^{\pm 1}, \ldots, t_{p-1}^{\pm 1}, t^{\pm 1}]$.
By Lemma \ref{lm:Torres}, we have
$$
{\mit \Delta}_{\widehat{\overline{L}}}(t_1, \ldots, t_{p-1}, t_1', \ldots, t_{q-1}', 1)
\doteq 
(t_1\cdots t_{p-1}-1){\mit \Delta}_{\widehat{L}}(t_1, \ldots, t_{p-1}, t_1', \ldots, t_{q-1}')
=0,
$$
and
$${\mit \Delta}_{\widehat{L}}(t_1, \ldots, t_{p-1}, t_1', \ldots, t_{q-1}')=0.
\quad
\qed$$

\begin{re}\label{re:MAlex}
{\rm
By Lemma \ref{lm:BrunnianAlex} (3), 
we obtain ${\mit \Delta}_{M_{\lambda}}(t_1, \ldots, t_{\lambda})=0$ for $\lambda \ge 4$
(cf.\ (\ref{eq:Alex})).}
\end{re}

\noindent
{\it Proof of Theorem \ref{th:Brunnian}.}\ 
It is easy to see that 
$\widehat{L}$ is a $(p+q-2)$-component link 
such that every component is an unknot.
We show that the Alexander polynomial of 
every proper sublink of $\widehat{L}$ with at least two components is zero.
We suppose that $i=p$ and $j=q$.
By Lemma \ref{lm:BrunnianAlex} (2) and (3), we have the result.
\qed

\bigskip

\begin{co}\label{co:cond}
Let $L=K_1\cup K_2$ be a $2$-component semi-Brunnian type link,
$K_2$ the characteristic component of $L$,
$L'=K_1'\cup \cdots \cup K_q'$ a $q$-component Brunnian type link with $q \ge 2$,
and $\widehat{L}=S_{2q}(L, L')$.
\begin{enumerate}
\item[(1)]
If there is another expression of $\widehat{L}$
which is a satellite of a $p'$-component semi-Brunnian type link 
with pattern a $q'$-component Brunnian type link such that
$p'\ge 3$ and $q'\ge 2$, then we have
$${\mit \Delta}_{\widehat{L}}(t_1, t_1', \ldots, t_{q-1}')=0.$$

\item[(2)]
If $|\mathrm{lk}(K_1, K_2)|\ge 2$ and
${\mit \Delta}_{L'}(t_1', \ldots, t_q')\ne 0$, then 
$\widehat{L}$ does not yield a lens space by any finite surgery.

\end{enumerate}
\end{co}

\noindent
{\it Proof.}\ 
(1) By Lemma \ref{lm:BrunnianAlex} (3), we have the result.

\medskip

\noindent
(2) By Lemma \ref{lm:BrunnianAlex} (2), we have
$${\mit \Delta}_{\widehat{L}}(t_1, t_1', \ldots, t_{q-1}')
\doteq
\left(\frac{t_1^k-1}{t_1-1}\right)^2
(t_1-1)(t_1'-1)\cdots (t_{q-1}'-1)f'(t_1', \ldots, t_{q-1}', t_1^k),$$
where $k=\mathrm{lk}(K_1, K_2)$.
We set
$$f(t_1, t_1', \ldots, t_{q-1}')
=\left(\frac{t_1^k-1}{t_1-1}\right)^2
f'(t_1', \ldots, t_{q-1}', t_1^k).$$
We prove that a Laurent polynomial which is obtained by
substituting $1$ to the variables in $f(t_1, t_1', \ldots, t_{q-1}')$
except one variable is not monic or
does not satisfy the conditions in Theorem \ref{th:OS} and Theorem \ref{th:Tange}.
It shows that $\widehat{L}$ does not yield a lens space by any finite surgery
by (\ref{eq:hatKi}), Theorem \ref{th:OS} and Theorem \ref{th:Tange}.

\medskip

If we substitute $1$ except $t_i'$ to $f(t_1, t_1', \ldots, t_{q-1}')$, 
then we obtain a Laurent polynomial of the form $k^2g_i(t_i')$
where $g_i(t_i')\in \mathbb{Z}[t_i', t_i'^{-1}]$.
If $|k|\ge 2$, then $f_i(t_i')$ is not monic.
If we substitute $1$ except $t_1$ to $f(t_1, t_1', \ldots, t_{q-1}')$, 
then we obtain a Laurent polynomial of the form 
\begin{equation}\label{eq:form}
\left(\frac{t_1^k-1}{t_1-1}\right)^2
f'(1, \ldots, 1, t_1^k).
\end{equation}
If $|k|\ge 2$, then the trace of the term $\left(\frac{t_1^k-1}{t_1-1}\right)^2$ is $-2$
and that of the term $f'(1, \ldots, 1, t_1^k)$ is $0$.
Since the degree of (\ref{eq:form}) is at least $2$,
both the trace of (\ref{eq:form}) and
that of ${\mit \Delta}_{\widehat{K}_1}(t_1)$ in (\ref{eq:hatKi}) are equal to $-2$.
It does not satisfy the condition in Theorem \ref{th:Tange}.
\qed

\bigskip

Finally, we pose the following questions :

\begin{qu}\label{qu:conj}
\begin{enumerate}
\item[(1)]
Are the Whitehead link and the Borromean rings the only non-trivial Brunnian links
and anannular Brunnian type links yielding lens spaces ?

\item[(2)]
Do every non-trivial Brunnian link and anannular Brunnian type link 
have Properties FP$'$, FP and FR ?

\end{enumerate}
\end{qu}

Corresponding to the cyclic surgery theorem \cite{CGLS},
we ask the following : 

\begin{qu}\label{qu:CGLS}
Are the cyclic surgery coefficients of the components of
a hyperbolic Brunnian link (resp.\ a hyperbolic anannular Brunnian type link)
integers except that of one component ?
\end{qu}

Theorem \ref{th:MT1}, Lemma \ref{lm:key2} (1) and Theorem \ref{th:MT4} are
partial answers for Question \ref{qu:CGLS} above.


\bigskip

{\noindent {\bf Acknowledgements}}\ 
The author would like to thank to
Yasuyoshi Tsutsumi for giving him useful advices
about the Casson-Walker invariant, 
Akio Kawauchi for pointing out relationship 
between the present result and 
his results \cite{Kw3} from the imitation theory \cite{Kw1, Kw2}, and
Motoo Tange for informing him about a result of M.~Hedden and L.~Watson \cite{HW}.

\bigskip

{\footnotesize
\par
\medskip
Teruhisa KADOKAMI\par
Department of Mathematics, East China Normal University,\par
Dongchuan-lu 500, Shanghai, 200241, China \par
{\tt mshj@math.ecnu.edu.cn} \par
{\tt kadokami2007@yahoo.co.jp} \par}

\end{document}